\documentclass[preprint,a4paper,12pt]{elsarticle}

\usepackage{lineno}

\usepackage{amsmath}
\usepackage{amsfonts}
\usepackage{amssymb}
\usepackage{amsthm}

\usepackage{graphicx}
\usepackage{subcaption}
\usepackage{multirow} % Table
\usepackage{float} % Position

\theoremstyle{definition}

\usepackage{xcolor}
\usepackage[colorlinks=true, allcolors=blue]{hyperref}

\newcommand{\bracketround}[1]{\left(#1\right)}

\newcommand{\abs}[1]{\left|#1\right|}

\begin{document}
	
\begin{frontmatter}
	
	\title{Do physics-informed neural networks (PINNs) need to be deep? Shallow PINNs using the Levenberg--Marquardt algorithm}
	
	\author[inst1,inst2]{Muhammad Luthfi Shahab\corref{mycorrespondingauthor}}
	\cortext[mycorrespondingauthor]{Corresponding author}
	\ead{shahab@its.ac.id}
	\author[inst1]{Imam Mukhlash}
	\author[inst2]{Hadi Susanto}
	
	\affiliation[inst1]{
		organization={Department of Mathematics, Institut Teknologi Sepuluh Nopember},
		% addressline={Address One}, 
		city={Surabaya},
		postcode={60111}, 
		% state={State One},
		country={Indonesia}
	}
	
	\affiliation[inst2]{
		organization={Department of Mathematics, Khalifa University of Science \& Technology},
		% addressline={Address One}, 
		city={Abu Dhabi},
		postcode={PO Box 127788}, 
		% state={State One},
		country={United Arab Emirates}
	}
	
	\begin{abstract}
		This work investigates shallow physics-informed neural networks (PINNs) for solving forward and inverse problems governed by nonlinear partial differential equations (PDEs). By formulating PINN training as a nonlinear least-squares problem, the Levenberg--Marquardt (LM) algorithm is used to efficiently optimize the network parameters.
		Exact analytical expressions for neural-network derivatives with respect to the input variables are derived, revealing the relationships between the network output and its spatial and temporal derivatives and providing a clearer interpretation of the PINN architecture. These expressions are then used to derive explicit formulas for the Jacobian matrix required by LM.
		The proposed approach is evaluated on the Burgers, Schr\"odinger, Allen--Cahn, and three-dimensional Bratu equations. Numerical results show that LM substantially outperforms BFGS, L-BFGS, and Adam in convergence speed, accuracy, and final loss values. Comparisons with deeper networks further demonstrate that shallow LM-PINNs can achieve higher accuracy with substantially fewer parameters, emphasizing the importance of considering network architecture and optimization strategy jointly. The explicit analytical Jacobian also provides computational and memory advantages that are particularly relevant to large-scale PINNs. Overall, these results suggest that, for a broad class of PDEs, shallow PINNs combined with effective second-order optimization can provide accurate and computationally efficient solutions to both forward and inverse problems.
		
	\end{abstract}
		
	
\end{frontmatter}

%\linenumbers

\section{Introduction}

In recent years, numerous deep learning-based approaches have been proposed for the numerical solution of partial differential equations (PDEs). Representative examples include physics-informed neural networks (PINNs) \cite{raissi2019physics}, deep learning methods for high-dimensional PDEs \cite{han2018solving, putri2024deep}, deep operator networks (DeepONet) \cite{lu2021learning}, and Fourier neural operators (FNOs) \cite{li2020fourier}. Among these approaches, PINNs have attracted particular attention because of their flexibility and ability to incorporate governing equations directly into the learning process.

Although the PINN framework was popularized by Raissi et al. \cite{raissi2019physics}, related ideas had been explored earlier, particularly for forward problems. Early contributions include the works of Lee and Kang \cite{lee1990neural}, Dissanayake and Phan-Thien \cite{dissanayake1994neural}, and Lagaris et al. \cite{lagaris1998artificial, lagaris2000neural}. Raissi et al. \cite{raissi2019physics} substantially extended these ideas to both forward and inverse problems for time-dependent PDEs using deep neural networks.
In PINNs, the governing PDEs and available data are incorporated into a composite loss function, which is typically minimized using first-order optimizers such as Adam \cite{kingma2014adam} or quasi-Newton methods such as L-BFGS \cite{liu1989limited}. More recently, BFGS, self-scaled BFGS (SSBFGS), and self-scaled Broyden (SSBroyden) methods have been shown to outperform L-BFGS in solution accuracy for PINNs \cite{kiyani2025optimizing}, motivating renewed interest in more effective optimization strategies.

Early PINN studies predominantly employed deep neural networks. However, deep architectures can incur substantial computational costs and training inefficiencies \cite{cuomo2022scientific}. To address these limitations, several alternatives based on shallow or parsimonious neural architectures have been proposed, including shallow neural networks \cite{hu2024shallow}, extreme learning machines (ELMs) \cite{huang2011extreme}, and random projection-based neural networks \cite{fabiani2023parsimonious}. Schiassi et al. \cite{schiassi2021extreme}, for example, introduced the extreme theory of functional connections (X-TFC), which employs a single-layer neural network trained using ELM. Dwivedi and Srinivasan \cite{dwivedi2020physics} proposed physics-informed extreme learning machines for linear PDEs, while Ramabathiran and Ramachandran \cite{ramabathiran2021spinn} developed sparse and partially interpretable PINNs (SPINNs) for ODEs and PDEs. Hu et al. \cite{hu2024shallow} constructed shallow PINNs for PDEs defined on static and evolving surfaces. Further related developments can be found in \cite{shi2024physics, rishi2024training, calabro2021extreme, fabiani2021numerical, fabiani2023parsimonious}.

Motivated by these developments, this work focuses on shallow PINNs by expressing the PINN residual equations as a nonlinear system and formulating the training process as a nonlinear least-squares problem. This formulation naturally enables the use of the Levenberg--Marquardt (LM) algorithm \cite{levenberg1944method, marquardt1963algorithm}, a classical second-order method for nonlinear least-squares problems. LM has been extensively studied for standard feed-forward neural networks. Hagan and Menhaj \cite{hagan1994training} demonstrated its effectiveness for training networks with a moderate number of parameters, while subsequent studies improved its memory efficiency and computational performance \cite{wilamowski2010improved, de2020stability}. LM has also been successfully applied in various engineering and control applications \cite{fu2014training, lv2017levenberg}.

Within the context of PINNs, applications of LM have so far been largely restricted to time-independent problems. Calandra et al. \cite{calandra2022multilevel} developed a multilevel LM strategy for stationary PDEs, while Yadav and Srinivasan \cite{yadavshallow} applied LM to singular perturbation problems. More recently, Shahab et al. \cite{shahab2024neural, shahab2025corrigendum, shahab2025physics} reformulated PINNs as nonlinear systems and employed LM to compute steady-state solutions and bifurcation diagrams, showing that shallow neural networks with only a few hundred parameters can accurately capture nonlinear dynamics. Applications of LM to time-dependent PDEs within the PINN framework remain comparatively limited. One notable example is Taylor et al. \cite{taylor2022optimizing}, who applied LM to a regression task and to a forward Burgers-equation problem.

This work develops a shallow-PINN framework based on LM optimization for nonlinear PDEs. The proposed approach combines shallow neural networks with exact analytical expressions for neural-network derivatives and explicit formulations of the Jacobian of the PDE residuals with respect to all trainable parameters. These expressions make explicit how the same network parameters determine both the solution and its spatial and temporal derivatives, thereby clarifying the parameter-sharing structure of PINNs. The resulting Jacobian is used directly in the LM iterations, and its relationship to the gradient of the PINN loss used by gradient-based and quasi-Newton methods is established. The framework is applied to both time-dependent and time-independent PDEs and to both forward and inverse problems. A central objective is to investigate whether deep neural networks are essential for high-accuracy PINNs or whether shallow architectures combined with an effective second-order optimization method can achieve comparable or superior performance.

The proposed approach is evaluated on four nonlinear PDEs: the Burgers, Schr\"odinger, Allen--Cahn \cite{raissi2019physics}, and three-dimensional Bratu equations \cite{shahab2025finite}. These problems provide complementary test cases covering time-dependent and time-independent settings, different nonlinearities, multiple spatial dimensions, and both real- and complex-valued solution representations. The Burgers and Schr\"odinger equations are considered as forward problems, whereas the Allen--Cahn and three-dimensional Bratu equations are treated as inverse problems, allowing both solution approximation and parameter identification to be assessed. 
In particular, the Schr\"odinger equation involves a complex-valued solution represented by multiple network outputs, while the three-dimensional Bratu equation provides a higher-dimensional spatial test problem. The performance of LM is compared with BFGS, L-BFGS, and Adam, including comparisons based on computational time to account for differences in per-iteration cost. Additional experiments examine robustness to observational noise and different initial guesses in the inverse problems. We further analyze the role of the LM damping parameter, computational and memory costs as the number of collocation points increases, and the relative contributions of Jacobian evaluation and the LM linear solve. Finally, comparisons with deeper L-BFGS networks provide additional evidence that high accuracy can be achieved with shallow architectures when paired with an effective second-order optimizer.

The remainder of this paper is organized as follows. Section 2 presents the PINN formulation as a nonlinear least-squares problem and describes the application of LM to forward and inverse problems. Section 3 reports the numerical results and performance comparisons. Section 4 concludes the paper and outlines directions for future research.

\section{Proposed PINNs}

This section presents the proposed PINN framework for forward and inverse problems and highlights its connections to and differences from the standard PINN approach.

\subsection{Shallow neural networks}

We consider shallow feed-forward neural networks with two hidden layers. For such architectures, most trainable parameters typically arise from the connections between the two hidden layers.

Throughout this paper, NN$(a,m_1,m_2,1)$ denotes a network with $a$ inputs, $m_1$ and $m_2$ neurons in the first and second hidden layers, respectively, and one output. Let
$W = \{W_1,B_1,W_2,B_2,W_3,B_3\}$
denote the set of all weights and biases. For a time-dependent PDE with one spatial variable, the network output $\tilde{u}(x,t,W)$ is defined as
\begin{equation} \label{eq_nn}
	\tilde{u}(x, t, W) = \sigma( \sigma( [x \; t] W_1 + B_1 )W_2 + B_2 )W_3 + B_3,
\end{equation}
where $\sigma(\cdot)$ is a nonlinear activation function and $[x \; t]$ is a row vector containing the spatial and temporal variables. The dimensions of the weights and biases are summarized in Tab.\ \ref{tab_weights}. For problems with additional spatial variables or without a temporal variable, only the input dimension, and hence the dimension of $W_1$, needs to be adjusted.

\begin{table}
	\small
	\centering
	\caption{Dimensions of the weight matrices and bias vectors for NN$(a, m_1, m_2, 1)$.}
	\label{tab_weights}
	\begin{tabular}{|c|c|l|}
		\hline
		Name  &       Size       & Description                                         \\ \hline
		$W_1$ &  $a \times m_1$  & Weights between the input and first hidden layers   \\
		$B_1$ &  $1 \times m_1$  & Biases of the first hidden layer                    \\
		$W_2$ & $m_1 \times m_2$ & Weights between the first and second hidden layers  \\
		$B_2$ &  $1 \times m_2$  & Biases of the second hidden layer                   \\
		$W_3$ &  $m_2 \times 1$  & Weights between the second hidden and output layers \\
		$B_3$ &   $1 \times 1$   & Bias of the output layer                            \\ \hline
	\end{tabular}
\end{table}

When neural networks are used to approximate PDE solutions, initial and boundary conditions can often be enforced directly through the solution representation. We therefore modify the network output as
\begin{equation} \label{eq_nn_enforce}
	u(x,t,W) = \tilde{u}(x,t,W) \cdot p(x,t) + q(x),
\end{equation}
where $p(x,t)$ and $q(x)$ are constructed so that $u(x,t,W)$ satisfies the prescribed initial and boundary conditions. Such exact enforcement can improve training efficiency and solution accuracy \cite{lu2021physics}. Here, $\tilde{u}$ and $u$ denote the network output before and after enforcing these conditions, respectively.
This formulation is particularly suitable when $p(x,t)$ and $q(x)$ can be constructed analytically, as is often possible for Dirichlet boundary conditions. Periodic conditions may require suitable transformations of the network inputs or outputs, while more complex geometries or Neumann conditions may require alternative boundary-enforcement strategies.

\subsection{PINNs for the forward problem}

Consider a general PDE of the form
\begin{equation} \label{eq_pde_forward}
	u_t + \mathcal{N}(u) = 0, \quad x \in \Omega, \quad t \in [0,T],
\end{equation}
with certain initial and boundary conditions, where $\mathcal{N}(\cdot)$ denotes a nonlinear differential operator. In the forward problem, the governing PDE is fully known, and the objective is to approximate the solution $u(x,t)$ over the spatiotemporal domain.

We approximate the solution by a neural network $u(x,t,W)$, where $W$ denotes the collection of all weights and biases. Let $\{(x_i,t_i)\}_{i=1}^n$ be $n$ collocation points sampled from the domain. Assuming that the initial and boundary conditions are enforced through the network construction \cite{lu2021physics}, the network parameters are typically determined by minimizing
\begin{equation} \label{eq_loss_forward}
	L(W) = \frac{1}{n} \sum_{i=1}^{n} (f_i(W))^2,
\end{equation}
where
\begin{equation} \label{eq_f_forward}
	f_i(W) = u_t(x_i, t_i, W) + \mathcal{N}(u(x_i, t_i, W)), \quad i=1,\dots,n,
\end{equation}
is the PDE residual at the $i$-th collocation point.

In the original PINN framework, this minimization problem is typically solved using the limited-memory BFGS (L-BFGS) algorithm \cite{liu1989limited}, while more recent studies have shown that BFGS can provide improved performance for PINN training \cite{kiyani2025optimizing}. In standard implementations, the spatial and temporal derivatives of $u(x,t,W)$ and the gradient of the loss with respect to the network parameters are computed using automatic differentiation \cite{baydin2018automatic}.

In this work, we instead assemble the PDE residuals into the nonlinear residual vector
\begin{equation} \label{eq_system_forward}
	F(W) = \begin{bmatrix}
		f_1(W) \\
		\vdots \\
		f_n(W)
	\end{bmatrix},
\end{equation}
and formulate PINN training as the nonlinear least-squares problem
\begin{equation} \label{eq_nls_forward}
	\min_W \frac{1}{2}\|F(W)\|_2^2.
\end{equation}
When an exact solution of the residual equations exists, this formulation corresponds to seeking $F(W)=\mathbf{0}$. We solve the resulting least-squares problem using the Levenberg--Marquardt (LM) algorithm \cite{levenberg1944method,marquardt1963algorithm}, which combines Gauss--Newton and gradient-descent behavior and is well suited to nonlinear least-squares problems.

Each LM iteration requires the Jacobian matrix $J$ of $F(W)$ with respect to all network parameters. Evaluating this Jacobian also requires the spatial and temporal derivatives of $u(x,t,W)$ that appear in the PDE residuals. Although these derivatives can be computed using automatic differentiation in deep learning frameworks such as TensorFlow \cite{abadi2016tensorflow}, LM is not readily supported in many modern frameworks. We therefore derive an exact analytical differentiation scheme for computing the required residual Jacobian. This formulation enables the direct application of LM to PINNs while also providing greater insight into their internal derivative structure and optimization process.

\subsubsection{Derivatives with respect to $x$ and $t$} \label{sec_derivatives_input}

In this subsection, we derive analytical expressions for the derivatives of the neural network functions $\tilde{u}$ and $u$ with respect to the input variables $x$ and $t$. We begin with the network output $\tilde{u}$ defined in Eq.\ \eqref{eq_nn}. For clarity, the network operations are written in layer-wise form:
\begin{equation} \label{eq_tilde_u}
	\begin{split}
		X & = [x \; t], \\
		h_1 & = X W_1 + B_1, \\
		H_1 & = \sigma(h_1), \\
		h_2 & = H_1 W_2 + B_2, \\
		H_2 & = \sigma(h_2), \\
		\tilde{u} & = H_2 W_3 + B_3.
	\end{split}
\end{equation}
Here, $h_\ell$ and $H_\ell$ denote the pre-activation and post-activation variables of the $\ell$-th hidden layer, respectively, and $\sigma(\cdot)$ is a nonlinear activation function, such as the hyperbolic tangent. In principle, any sufficiently smooth activation function can be used, provided that the derivatives required by the PDE exist. In particular, $\sigma\in C^\infty$ guarantees derivatives of arbitrary order.

We first compute the partial derivative of $\tilde{u}$ with respect to $t$. Since $X=[x\;t]$, its derivative with respect to $t$ is given by $X_t=[0\;1]$, which leads to
\begin{equation} \label{eq_tilde_ut}
	\begin{split}
		X_t & = [0 \; 1], \\
		h_{1,t} & = X_t W_1, \\
		H_{1,t} & = \sigma'(h_1) \cdot h_{1,t}, \\
		h_{2,t} & = H_{1,t} W_2, \\
		H_{2,t} & = \sigma'(h_2) \cdot h_{2,t}, \\
		\tilde{u}_t & = H_{2,t} W_3. 
	\end{split}
\end{equation}
Similarly, using $X_x=[1\;0]$, the derivative of $\tilde{u}$ with respect to $x$ is
\begin{equation} \label{eq_tilde_ux}
	\begin{split}
		X_x & = [1 \; 0], \\
		h_{1,x} & = X_x W_1, \\
		H_{1,x} & = \sigma'(h_1) \cdot h_{1,x}, \\
		h_{2,x} & = H_{1,x} W_2, \\
		H_{2,x} & = \sigma'(h_2) \cdot h_{2,x}, \\
		\tilde{u}_x & = H_{2,x} W_3.
	\end{split}
\end{equation}
The second derivative with respect to $x$ is then
\begin{equation} \label{eq_tilde_uxx}
	\begin{split}
		X_{xx} & = [0 \; 0], \\
		h_{1,xx} & = X_{xx} W_1, \\
		H_{1,xx} & = \sigma''(h_1) \cdot h_{1,x}^2 + \sigma'(h_1) \cdot h_{1,xx}, \\
		h_{2,xx} & = H_{1,xx} W_2, \\
		H_{2,xx} & = \sigma''(h_2) \cdot h_{2,x}^2 + \sigma'(h_2) \cdot h_{2,xx}, \\
		\tilde{u}_{xx} & = H_{2,xx} W_3.
	\end{split}
\end{equation}
Although $X_{xx}$ and $h_{1,xx}$ vanish for the present input representation, we retain these terms for generality. If the inputs are transformed nonlinearly, for example through $\sin(x)$ or $\cos(x)$, their second derivatives need not vanish. The same formulation therefore extends naturally to such cases.

When the initial and boundary conditions are enforced through the transformation in Eq.\ \eqref{eq_nn_enforce}, the resulting approximation and its derivatives are
\begin{equation} \label{eq_u_ut_ux_uxx}
	\begin{split}
		u & = \tilde{u} \cdot p(x,t) + q(x), \\
		u_t & = \tilde{u}_t \cdot p(x,t) + \tilde{u} \cdot p_t(x,t), \\
		u_x & = \tilde{u}_x \cdot p(x,t) + \tilde{u} \cdot p_x(x,t) + q_x(x), \\
		u_{xx} & = \tilde{u}_{xx} \cdot p(x,t) + 2 \tilde{u}_x \cdot p_x(x,t) + \tilde{u} \cdot p_{xx}(x,t) + q_{xx}(x).
	\end{split}
\end{equation}
For brevity, only derivatives up to second order in $x$ are presented, although higher-order derivatives can be derived analogously. The procedure also extends directly to deeper network architectures. Once $u$, $u_t$, $u_x$, and $u_{xx}$ have been computed, they can be substituted into the residual vector $F(W)$ defined in Eq.\ \eqref{eq_system_forward}.

The same computations extend naturally to a vectorized implementation. Let $\mathbf{x}$ and $\mathbf{t}$ denote vectors containing $n$ collocation points. Then
\begin{equation}
	X = [\mathbf{x} \; \mathbf{t}], \quad
	X_t = [\mathbf{0} \; \mathbf{1}], \quad
	X_x = [\mathbf{1} \; \mathbf{0}], \quad
	X_{xx} = [\mathbf{0} \; \mathbf{0}],
\end{equation}
where $\mathbf{0}$ and $\mathbf{1}$ are vectors of zeros and ones of length $n$. Hence, $X$, $X_t$, $X_x$, and $X_{xx}$ are all $n\times2$ matrices. For a network with architecture NN$(2,m_1,m_2,1)$, $h_1$ and $H_1$ have dimensions $n\times m_1$, whereas $h_2$ and $H_2$ have dimensions $n\times m_2$. In Eqs.\ \eqref{eq_tilde_u}--\eqref{eq_u_ut_ux_uxx}, the operator $\cdot$ denotes element-wise multiplication, while products such as $XW_1$ denote standard matrix multiplication. Terms such as $h_{1,x}^2$ denote element-wise squaring.

\begin{figure}
	\centering
	\includegraphics[width=\textwidth]{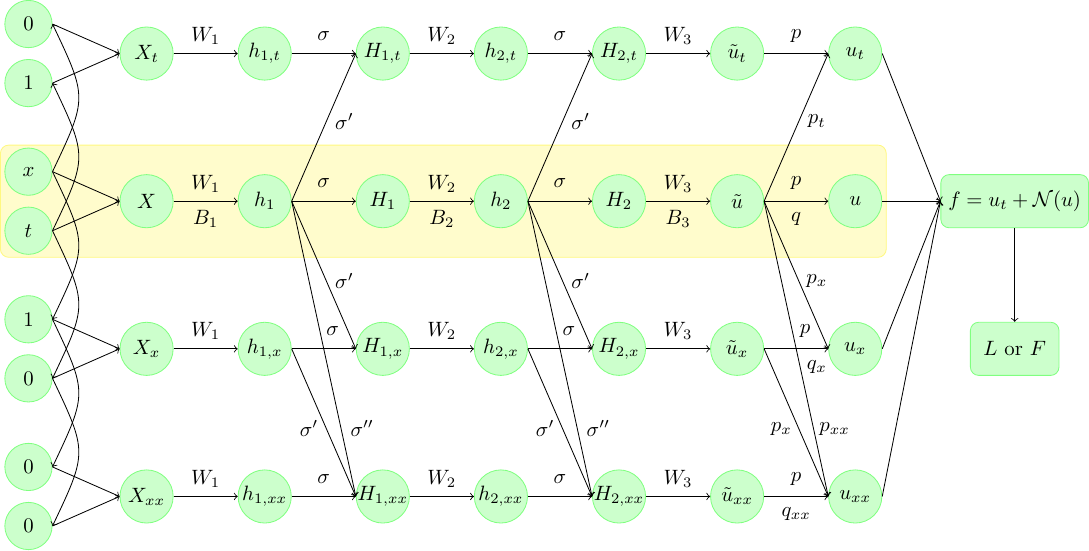}
	\caption{Explicit representation of the internal PINN architecture, showing how the core network $\tilde{u}$ (yellow) and its spatial and temporal derivatives are constructed analytically from shared weights and biases.}
	\label{fig_pinns}
\end{figure}

Fig.\ \ref{fig_pinns} illustrates the relationships between $\tilde{u}$ and its derivatives $\tilde{u}_t$, $\tilde{u}_x$, and $\tilde{u}_{xx}$. 
The yellow block represents the forward pass used to compute $\tilde{u}$, while the derivative calculations reuse the same network weights and biases.
For clarity, each hidden layer is represented by a single node, allowing the figure to emphasize the data flow and parameter-dependency structure rather than individual neurons.

To the best of our knowledge, an explicit illustration of the internal PINN architecture and its analytical connections to spatial and temporal derivatives has not been presented in this form in previous studies. Most existing approaches rely on automatic differentiation, which can obscure the dependency structure between the network parameters and the resulting differential operators. By deriving and visualizing these relationships explicitly, the proposed formulation provides additional insight into the internal mechanics of PINNs and clarifies how the solution and its derivatives are constructed from shared network parameters.

\subsubsection{Jacobian}

As discussed earlier, the PINN residuals are assembled into a nonlinear residual vector, and the resulting nonlinear least-squares problem is solved using LM. Therefore, computation of the Jacobian matrix is a key step in each LM iteration.

To construct the Jacobian, all neural network weights and biases, namely $W_1$, $B_1$, $W_2$, $B_2$, $W_3$, and $B_3$, are reshaped into row vectors and concatenated into a single column vector
\begin{equation} \label{eq_weight_ordered}
	W = \begin{bmatrix}
		W_1 & B_1 & W_2 & B_2 & W_3 & B_3
	\end{bmatrix}^T.
\end{equation}
Let $m$ denote the total number of trainable parameters, so that $W \in \mathbb{R}^{m}$. For $n$ collocation points $\{(x_i,t_i)\}_{i=1}^n$, the Jacobian matrix $J \in \mathbb{R}^{n \times m}$ is defined as
\begin{equation}
	J = \begin{bmatrix}
		\frac{\partial f_1}{\partial w_1} & \cdots & \frac{\partial f_1}{\partial w_m} \\
		\vdots & \ddots & \vdots \\
		\frac{\partial f_n}{\partial w_1} & \cdots & \frac{\partial f_n}{\partial w_m}
	\end{bmatrix},
\end{equation}
where $w_r$, $r=1,\dots,m$, denotes the $r$-th component of $W$.
The $(i,r)$-th entry of the Jacobian is given by
\begin{equation} \label{eq_jacobian}
	\begin{split}
		J_{i,r} & = \frac{\partial f_i}{\partial w_r} \\
		& = \frac{\partial ( u_t(x_i, t_i, W) + \mathcal{N}(u(x_i, t_i, W)) )}{\partial w_r} \\
		& = \frac{\partial u_t(x_i, t_i, W)}{\partial w_r} + \mathcal{N}_{w_r}(u(x_i, t_i, W)),
	\end{split}
\end{equation}
where $\mathcal{N}_{w_r}(\cdot)$ denotes the partial derivative of the nonlinear operator $\mathcal{N}$ with respect to $w_r$. This derivative generally involves derivatives of $u$, $u_x$, and $u_{xx}$ with respect to $w_r$.

As an illustrative example, consider
\begin{equation}
	\mathcal{N}(u) = u u_x - \frac{0.01}{\pi} u_{xx}.
\end{equation}
Its derivative with respect to $w_r$ is
\begin{equation}
	\mathcal{N}_{w_r}(u) = \frac{\partial u}{\partial w_r} u_x 
	+ u \frac{\partial u_x}{\partial w_r}
	- \frac{0.01}{\pi} \frac{\partial u_{xx}}{\partial w_r}.
\end{equation}
Thus, constructing the Jacobian requires explicit expressions for $\frac{\partial u}{\partial w_r}$, $\frac{\partial u_t}{\partial w_r}$, $\frac{\partial u_x}{\partial w_r}$, and $\frac{\partial u_{xx}}{\partial w_r}$.

It is worth noting that these derivatives depend on the position of $w_r$ within the network. Parameters closer to the output layer generally yield simpler expressions, whereas those near the input layer require more extensive chain-rule expansions. This also illustrates a computational advantage of shallow architectures, since for a fixed number of parameters, shallower networks generally require fewer operations to evaluate the Jacobian. \ref{sec_derivatives_weight} provides explicit formulas for the derivatives of $\tilde{u}$, $\tilde{u}_t$, $\tilde{u}_x$, and $\tilde{u}_{xx}$ with respect to $w_r$ for each parameter location in the network.

Since the initial and boundary conditions are enforced through Eq.\ \eqref{eq_u_ut_ux_uxx}, the corresponding derivatives of $u$ and its temporal and spatial derivatives are
\begin{equation} \label{eq_u_ut_ux_uxx_wr}
	\begin{split}
		\frac{\partial u}{\partial w_r} 
		& = \frac{\partial \tilde{u}}{\partial w_r} \cdot p(x,t), \\
		\frac{\partial u_t}{\partial w_r} 
		& = \frac{\partial \tilde{u}_t}{\partial w_r} \cdot p(x,t) 
		+ \frac{\partial \tilde{u}}{\partial w_r} \cdot p_t(x,t), \\
		\frac{\partial u_x}{\partial w_r} 
		& = \frac{\partial \tilde{u}_x}{\partial w_r} \cdot p(x,t) 
		+ \frac{\partial \tilde{u}}{\partial w_r} \cdot p_x(x,t), \\
		\frac{\partial u_{xx}}{\partial w_r} 
		& = \frac{\partial \tilde{u}_{xx}}{\partial w_r} \cdot p(x,t) 
		+ 2 \frac{\partial \tilde{u}_x}{\partial w_r} \cdot p_x(x,t) 
		+ \frac{\partial \tilde{u}}{\partial w_r} \cdot p_{xx}(x,t).
	\end{split}
\end{equation}
These derivatives are substituted into Eq.\ \eqref{eq_jacobian} to assemble the full Jacobian matrix $J$. The LM parameter correction is then computed as
\begin{equation} \label{eq_LM_update}
	\Delta W = (J^T J + \mu I)^{-1} J^T F,
\end{equation}
where $\mu>0$ is the damping parameter that controls the transition between Gauss--Newton and gradient-descent behavior \cite{hagan1994training}. The parameter vector is then updated according to
\begin{equation}
	W_{\mathrm{new}} = W - \Delta W.
\end{equation}

The damping parameter $\mu$ is adjusted adaptively by a factor of 10. If a trial step reduces the sum of squared errors (SSE), the step is accepted and $\mu$ is decreased according to $\mu\leftarrow\mu/10$. Otherwise, the step is rejected and $\mu$ is increased according to $\mu\leftarrow10\mu$. This process is repeated until an improving trial step is obtained, after which the accepted parameters are used in the next iteration.

\subsubsection{Connection between Jacobian and gradient}

In standard PINNs, when the loss function is defined as in Eq.\ \eqref{eq_loss_forward} and optimized using quasi-Newton methods such as BFGS or L-BFGS, the gradient of the loss with respect to the network parameters must be evaluated at each iteration.

Let the weights and biases be reshaped and ordered as in Eq.\ \eqref{eq_weight_ordered}. Recall from Eq.\ \eqref{eq_jacobian} that
$J_{i,r} = \frac{\partial f_i}{\partial w_r}$,
where $f_i$ denotes the residual at the $i$th collocation point and $w_r$ is the $r$-th component of the parameter vector $W$. The $r$-th component of the gradient vector $G$ is therefore
\begin{equation}
	\frac{\partial L}{\partial w_r} = \frac{2}{n} \sum_{i=1}^{n} \frac{\partial f_i}{\partial w_r} f_i = \frac{2}{n} \sum_{i=1}^{n} J_{i,r} f_i.
\end{equation}
In matrix form, this can be written compactly as
\begin{equation}
	G = \frac{2}{n} J^T F.
\end{equation}
where $J$ is the Jacobian matrix and $F$ is the residual vector defined in Eq.\ \eqref{eq_system_forward}.

This relation shows that the residual Jacobian required by LM is directly related to the gradient used by BFGS and other gradient-based optimization methods. In standard automatic-differentiation implementations, the product $J^T F$ can be evaluated without explicitly forming the full Jacobian. LM, however, requires the Jacobian more explicitly to construct the approximate Hessian $J^T J$. Thus, assuming comparable costs for residual and derivative evaluations, an important difference between LM and BFGS lies in the additional linear-algebra operations required by their respective parameter updates.

\subsection{PINNs for the inverse problem}

We next consider inverse problems. Since the formulation closely follows that of the forward problem, only the main differences are described. Consider a general PDE of the form
\begin{equation} \label{eq_pde_inverse}
	u_t + \mathcal{N}(u, \lambda) = 0, \quad x \in \Omega, \quad t \in [0,T],
\end{equation}
with certain initial and boundary conditions,
where $\mathcal{N}(\cdot,\lambda)$ denotes a nonlinear differential operator and $\lambda=\{\lambda_1,\lambda_2,\dots\}$ represents the unknown PDE parameters. Given partial observations of the solution and the known structure of the PDE, the objective is to identify $\lambda$ while simultaneously approximating $u(x,t)$.

The solution is approximated by a neural network $u(x,t,W)$ with trainable weights and biases $W$. Let $\{(x_i,t_i),u_i^*\}_{i=1}^n$ denote $n$ collocation points and their corresponding observed values. Assuming that the initial and boundary conditions are incorporated into the network architecture \cite{lu2021physics}, the network parameters $W$ and the unknown PDE parameters $\lambda$ can be determined by minimizing
\begin{equation} \label{eq_loss_inverse}
	L(W, \lambda) = \frac{1}{n} \sum_{i=1}^{n} (f_i(W, \lambda))^2 + \frac{1}{n} \sum_{i=1}^{n} (g_i(W))^2 ,
\end{equation}
where
\begin{equation} \label{eq_f_inverse}
	f_i(W, \lambda) = u_t(x_i, t_i, W) + \mathcal{N}(u(x_i, t_i, W), \lambda),
	\quad i=1,\dots,n,
\end{equation}
is the PDE residual, and
\begin{equation} \label{eq_g_inverse}
	g_i(W) =  u(x_i, t_i, W) - u^*_i, \quad i=1,\dots,n.
\end{equation}
measures the discrepancy between the predicted and observed solutions.

As in the forward problem, the residuals are assembled into
\begin{equation} \label{eq_system_inverse}
	F(W, \lambda) = \begin{bmatrix}
		F_1(W, \lambda) \\
		F_2(W)
	\end{bmatrix},
\end{equation}
where
\begin{equation} \label{eq_system_inverse1}
	F_1(W, \lambda) = \begin{bmatrix}
		f_1(W, \lambda) \\
		\vdots \\
		f_n(W, \lambda)
	\end{bmatrix} , \quad
	F_2(W) = \begin{bmatrix}
		g_1(W) \\
		\vdots \\
		g_n(W)
	\end{bmatrix} .
\end{equation}
The corresponding nonlinear least-squares problem is then solved using LM \cite{levenberg1944method,marquardt1963algorithm}.

The derivatives of the neural network with respect to $x$ and $t$ are computed in the same manner as in the forward problem. To include the unknown PDE parameters, we define the augmented parameter vector
\begin{equation} \label{eq_weight_ordered_lambda}
	W \cup \lambda  = \begin{bmatrix}
		W_1 & B_1 & W_2 & B_2 & W_3 & B_3 & \lambda
	\end{bmatrix}^T.
\end{equation}
If $m$ denotes the number of network weights and biases and $m_\lambda$ the number of unknown PDE parameters, then $W \cup \lambda\in\mathbb{R}^{m+m_\lambda}$.

The main distinction from the forward problem lies in the structure of the Jacobian. For $n$ collocation points, $J\in\mathbb{R}^{2n\times(m+m_\lambda)}$ has the block form
\begin{equation}
	J = \begin{bmatrix}
		J_1 & J_3 \\
		J_2 & \textbf{0}
	\end{bmatrix},
\end{equation}
where
\begin{equation}
	J_1 = \begin{bmatrix}
		\frac{\partial f_1}{\partial w_1} & \cdots & \frac{\partial f_1}{\partial w_m}  \\
		\vdots & \ddots & \vdots \\
		\frac{\partial f_n}{\partial w_1} & \cdots & \frac{\partial f_n}{\partial w_m} 
	\end{bmatrix},
	J_2 = \begin{bmatrix}
		\frac{\partial g_1}{\partial w_1} & \cdots & \frac{\partial g_1}{\partial w_m}  \\
		\vdots & \ddots & \vdots \\
		\frac{\partial g_n}{\partial w_1} & \cdots & \frac{\partial g_n}{\partial w_m} 
	\end{bmatrix},
	J_3 = \begin{bmatrix}
		\frac{\partial f_1}{\partial \lambda_1} & \cdots & \frac{\partial f_1}{\partial \lambda_{m_\lambda}}  \\
		\vdots & \ddots & \vdots \\
		\frac{\partial f_n}{\partial \lambda_1} & \cdots & \frac{\partial f_n}{\partial \lambda_{m_\lambda}} 
	\end{bmatrix}.
\end{equation}

For $r=1,\dots,m$, the derivatives with respect to the network parameters are
\begin{equation}
	\begin{split}
		\frac{\partial f_i}{\partial w_r} 
		& = \frac{\partial ( u_t(x_i, t_i, W) + \mathcal{N}(u(x_i, t_i, W), \lambda) )}{\partial w_r} \\
		& = \frac{\partial u_t(x_i, t_i, W)}{\partial w_r} + \mathcal{N}_{w_r}(u(x_i, t_i, W), \lambda),
	\end{split}
\end{equation}
and,
\begin{equation}
	\begin{split}
		\frac{\partial g_i}{\partial w_r} 
		& = \frac{\partial ( u(x_i, t_i, W) - u^*_i )}{\partial w_r} \\
		& = \frac{\partial u(x_i, t_i, W)}{\partial w_r}.
	\end{split}
\end{equation}
For $\lambda_r$, $r=1,\dots,m_\lambda$,
\begin{equation}
	\begin{split}
		\frac{\partial f_i}{\partial \lambda_r} 
		& = \frac{\partial ( u_t(x_i, t_i, W) + \mathcal{N}(u(x_i, t_i, W), \lambda) )}{\partial \lambda_r} \\
		& = \mathcal{N}_{\lambda_r}(u(x_i, t_i, W)).
	\end{split}
\end{equation}
For example, if
\begin{equation}
	\mathcal{N}(u) = \lambda_1 u u_x + \lambda_2 u_{xx},
\end{equation}
then
\begin{equation}
	\mathcal{N}_{\lambda_1}(u) = u u_x, 
	\quad
	\mathcal{N}_{\lambda_2}(u) = u_{xx}.
\end{equation}
The derivatives of $u$, $u_t$, $u_x$, and $u_{xx}$ with respect to $w_r$ are computed using the same analytical differentiation procedure described for the forward problem.

Finally, the gradient of the loss function \eqref{eq_loss_inverse} can also be expressed in terms of the Jacobian. For $r=1,\dots,m$,
\begin{equation}
	\frac{\partial L}{\partial w_r} 
	= \frac{2}{n} \sum_{i=1}^{n} \frac{\partial f_i}{\partial w_r} f_i 
	+ \frac{2}{n} \sum_{i=1}^{n} \frac{\partial g_i}{\partial w_r} g_i.
\end{equation}
while for $r=1,\dots,m_\lambda$,
\begin{equation}
	\frac{\partial L}{\partial \lambda_r} 
	= \frac{2}{n} \sum_{i=1}^{n} \frac{\partial f_i}{\partial \lambda_r} f_i.
\end{equation}
In compact form,
\begin{equation}
	G = \begin{bmatrix}
		\frac{2}{n} J_1^T F_1 + \frac{2}{n} J_2^T F_2 \\
		\frac{2}{n} J_3^T F_1
	\end{bmatrix}.
\end{equation}

\section{Experimental results}

We compare the performance of LM, BFGS, L-BFGS, and Adam for training PINNs. For BFGS, L-BFGS, and Adam, the gradients are computed using the Jacobian--gradient relationships derived in the previous section. Because an LM iteration requires solving a linear system and is generally more expensive than an iteration of the other methods, the maximum number of iterations is set to 4000 for LM and 20000 for BFGS, L-BFGS, and Adam. These limits provide sufficient optimization steps while keeping the overall computational costs within a comparable range and allowing each method to approach a stable optimization regime.

For LM, we use MATLAB's built-in \texttt{fsolve} function with the \texttt{'Algorithm'} option set to \texttt{'levenberg-marquardt'}. BFGS and L-BFGS are implemented using MATLAB's \texttt{fminunc} function with the \texttt{'HessianApproximation'} option set to \texttt{'bfgs'} and \texttt{\{'lbfgs', 100\}}, respectively, while Adam is implemented manually. Additional optimization settings are provided in \ref{sec_parameters}. Since LM in \texttt{fsolve} is applied to the nonlinear residual vector, whereas the other methods minimize a scalar-valued objective function, the corresponding PINN loss is evaluated after each LM iteration using Eq.\ \eqref{eq_loss_forward} for forward problems and Eq.\ \eqref{eq_loss_inverse} for inverse problems. This enables all methods to be compared using the same loss-based metric.

In all experiments, the hyperbolic tangent (\texttt{tanh}) activation function is used, as is common in PINNs. The same collocation points are also used for all optimization methods so that differences in performance can be attributed to the optimization algorithms rather than to variations in the training data.
All simulations are performed in MATLAB R2024a on a computer equipped with a six-core Intel Core i5-10400 CPU.

\subsection{Forward problem: Burgers equation}

We first consider the forward problem for the one-dimensional Burgers equation studied by Raissi et al. \cite{raissi2019physics},
\begin{equation} \label{eq_burgers1}
	\begin{split}
		& u_t + u u_x - \frac{0.01}{\pi} u_{xx} = 0, \quad x \in [-1,1], \quad t \in [0,1], \\
		& u(x, 0) = - \sin(\pi x), \\
		& u(-1, t) = u(1, t) = 0.
	\end{split}
\end{equation}
To enforce the initial and boundary conditions, we use the trial solution
\begin{equation} \label{eq_burgers_enforce}
	u = \tilde{u} (1 - x^2) t - \sin(\pi x).
\end{equation}
which automatically satisfies all prescribed conditions. For training, $n=10000$ collocation points are sampled from the spatiotemporal domain $[-1,1]\times[0,1]$. We initially use the shallow architecture NN$(2,25,25,1)$, which contains 751 trainable weights and biases.

\begin{figure}
	\centering
	
	\begin{subfigure}[b]{\textwidth}
		\centering
		\includegraphics[width=0.9\textwidth]{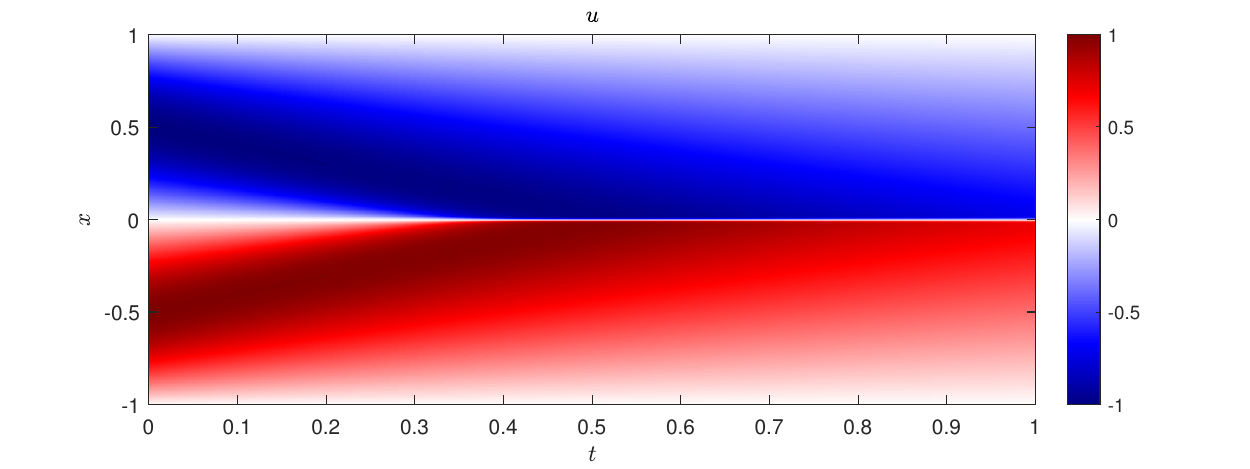}
		\caption{}
		\label{fig_forward_burgers_solution1}
	\end{subfigure}
	
	\begin{subfigure}[b]{\textwidth}
		\centering
		\includegraphics[width=\textwidth]{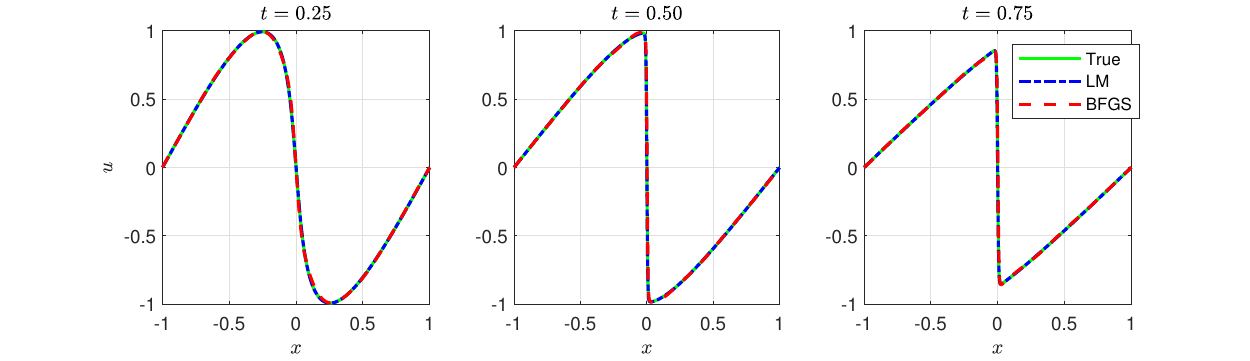}
		\caption{}
		\label{fig_forward_burgers_solution2}
	\end{subfigure}
	
	\begin{subfigure}[b]{\textwidth}
		\centering
		\includegraphics[width=\textwidth]{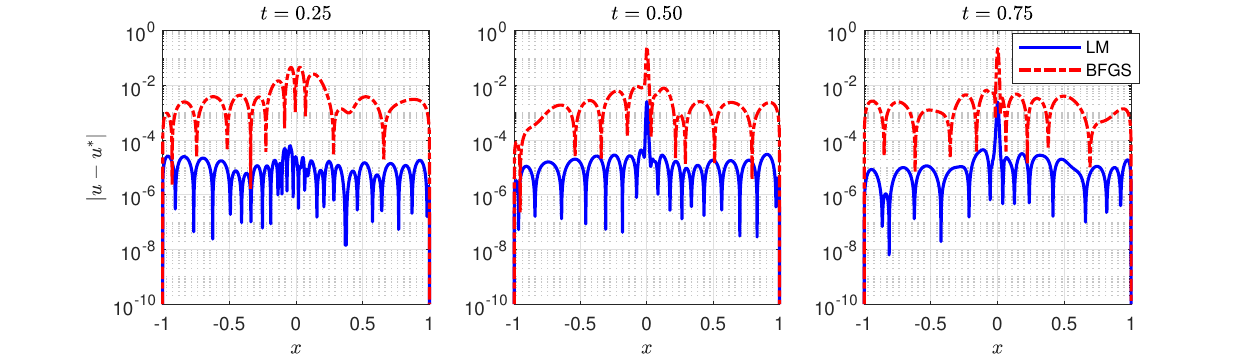}
		\caption{}
		\label{fig_forward_burgers_solution3}
	\end{subfigure}
	
	\begin{subfigure}[b]{0.49\textwidth}
		\centering
		\includegraphics[width=0.9\textwidth]{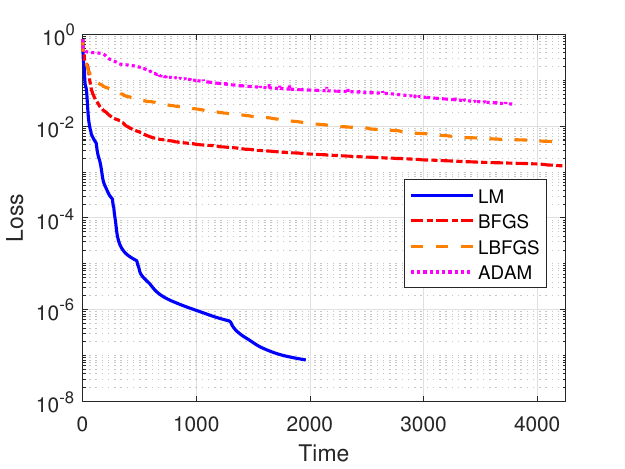}
		\caption{}
		\label{fig_forward_burgers_loss_time}
	\end{subfigure}
	
	\caption{(a) Neural network solution of the Burgers equation obtained using LM.
		(b) Solutions computed by LM and BFGS at selected time slices.
		(c) Absolute errors corresponding to the LM and BFGS solutions.
		(d) Training loss histories of LM, BFGS, L-BFGS, and Adam.}
	\label{fig_forward_burgers}
\end{figure}

Fig.\ \ref{fig_forward_burgers_solution1} shows the spatiotemporal solution obtained using LM. In Fig.\ \ref{fig_forward_burgers_solution2}, the LM and BFGS solutions at $t=0.25$, $0.50$, and $0.75$ are compared with the reference solution from \cite{basdevant1986spectral}. The corresponding absolute errors are shown in Fig.\ \ref{fig_forward_burgers_solution3}, where LM achieves errors approximately three orders of magnitude smaller than those of BFGS. Fig.\ \ref{fig_forward_burgers_loss_time} shows the evolution of the loss with computational time. LM reaches a substantially lower final loss of $8.1\times10^{-8}$ while requiring only about half the computational time of BFGS, L-BFGS, and Adam.
To further quantify the solution accuracy, we compute the relative $L_2$ error against the reference solution on a uniform $1001\times1001$ grid with $\Delta x=2\times10^{-3}$ and $\Delta t=10^{-3}$. The resulting errors are $2.0\times10^{-4}$ for LM and $2.2\times10^{-2}$ for BFGS, further confirming the higher accuracy of LM for this problem.

\begin{table}
	\small
	\centering
	\caption{Comparison of loss values, relative $L_2$ errors, and computational times for different network sizes trained using LM for the Burgers equation.}
	\label{tab_forward_burgers_1}
	\begin{tabular}{|l|c|c|c|}
		\hline
		Network         &       Loss        & Relative $L_2$ error & Time \\ \hline
		NN$(2,5,5,1)$   & $2.2\cdot10^{-2}$ &  $5.3\cdot10^{-2}$   & 108  \\ \hline
		NN$(2,10,10,1)$ & $1.0\cdot10^{-3}$ &  $2.4\cdot10^{-2}$   & 262  \\ \hline
		NN$(2,15,15,1)$ & $1.6\cdot10^{-5}$ &  $2.3\cdot10^{-3}$   & 750  \\ \hline
		NN$(2,20,20,1)$ & $1.2\cdot10^{-6}$ &  $5.4\cdot10^{-4}$   & 1351 \\ \hline
		NN$(2,25,25,1)$ & $8.1\cdot10^{-8}$ &  $2.0\cdot10^{-4}$   & 1962 \\ \hline
	\end{tabular}
\end{table}

Tab.\ \ref{tab_forward_burgers_1} summarizes the loss values, relative $L_2$ errors, and computational times for several network sizes trained using LM. As expected, increasing the number of neurons improves accuracy. Notably, even NN$(2,15,15,1)$ achieves good accuracy at substantially lower computational cost, indicating that moderately sized shallow networks can provide an effective balance between accuracy and efficiency.

\subsection{Forward problem: nonlinear Schr\"odinger equation}

We next consider the nonlinear Schr\"odinger equation \cite{raissi2019physics},
\begin{equation} \label{eq_nls1}
	\begin{split}
		& i u_t + \frac{1}{2} u_{xx} + \abs{u}^2 u = 0, \quad x \in [-5,5], \quad t \in [0, \pi/2], \\
		& u(x, 0) = 2 \; \text{sech} (x) \cos^2 \bracketround{\frac{\pi}{10} x}, \\
		& u(-5, t) = u(5, t), \\
		& u_x(-5, t) = u_x(5, t).
	\end{split}
\end{equation}
where $u$ is complex-valued.

We decompose $u$ into its real and imaginary parts,
\begin{equation}
	u(x,t) = v(x,t) + i w(x,t),
\end{equation}
which yields the real-valued system
\begin{equation} \label{eq_nls2}
	\begin{split}
		& - w_t + \frac{1}{2} v_{xx} + (v^2 + w^2) v = 0, \\
		& v_t + \frac{1}{2} w_{xx} + (v^2 + w^2) w = 0.
	\end{split}
\end{equation}
Unlike \cite{raissi2019physics}, which used a single PINN with two outputs, we employ two separate PINNs to approximate $\tilde{v}$ and $\tilde{w}$ and train them simultaneously. The initial conditions are enforced through
\begin{equation}
	v = \tilde{v} t + 2 \; \text{sech} (x) \cos^2 \bracketround{\frac{\pi}{10} x},
	\quad
	w = \tilde{w} t.
\end{equation}
Periodic boundary conditions are imposed by mapping the spatial coordinate to $\sin(2\pi x/P)$ and $\cos(2\pi x/P)$ with $P=10$ \cite{lu2021physics}. The resulting network input is
\begin{equation}
	X  = [\sin(2\pi x/P) \; \cos(2\pi x/P) \; t].
\end{equation}
We again use $n=10000$ collocation points and begin with two NN$(3,25,25,1)$ networks, one for $v$ and one for $w$, resulting in 1552 total trainable parameters.

\begin{figure}
	\centering
	
	\begin{subfigure}[b]{\textwidth}
		\centering
		\includegraphics[width=0.9\textwidth]{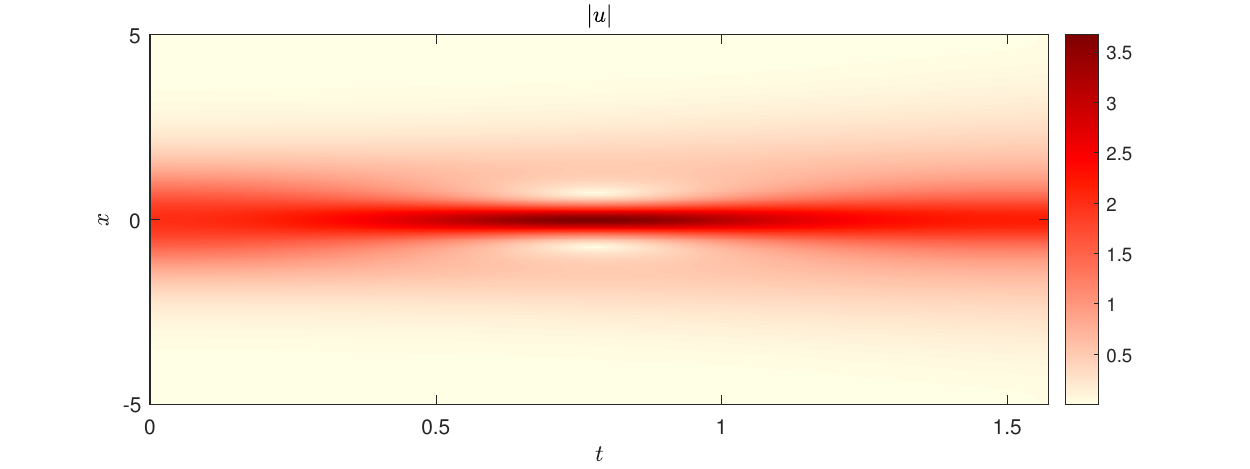}
		\caption{}
		\label{fig_forward_nls_solution1}
	\end{subfigure}
	
	\begin{subfigure}[b]{\textwidth}
		\centering
		\includegraphics[width=\textwidth]{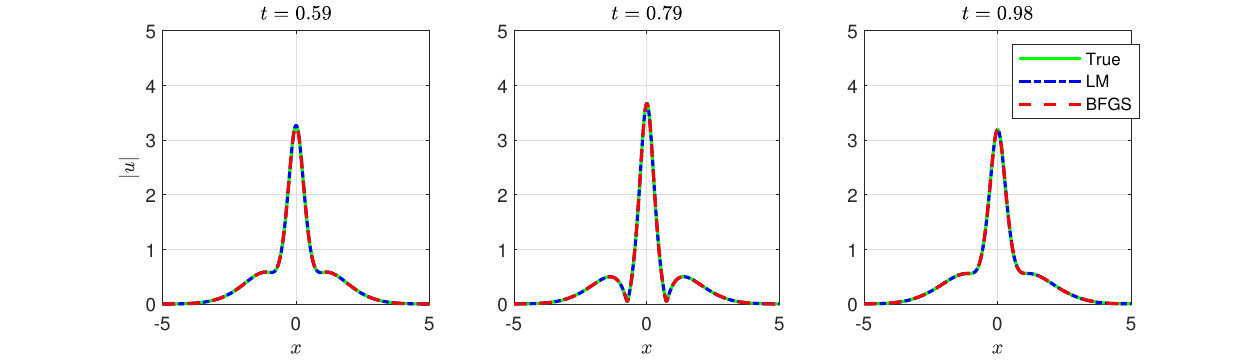}
		\caption{}
		\label{fig_forward_nls_solution2}
	\end{subfigure}
	
	\begin{subfigure}[b]{\textwidth}
		\centering
		\includegraphics[width=\textwidth]{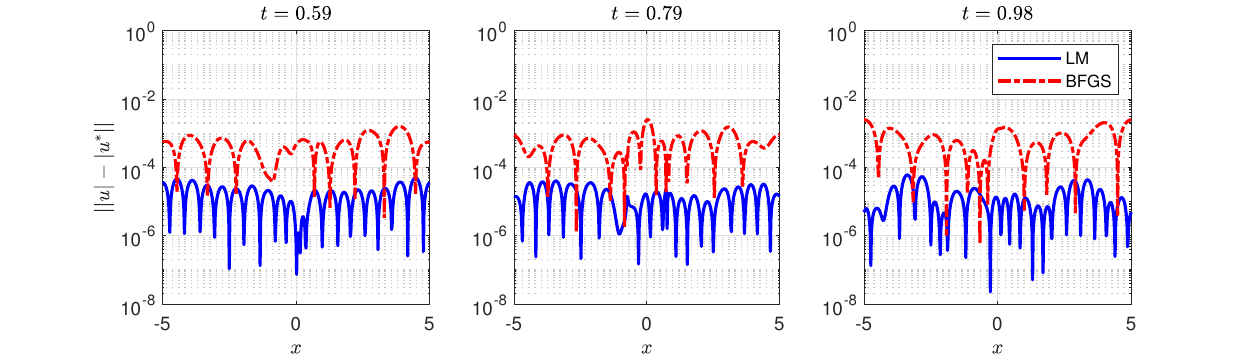}
		\caption{}
		\label{fig_forward_nls_solution3}
	\end{subfigure}
	
	\begin{subfigure}[b]{0.49\textwidth}
		\centering
		\includegraphics[width=0.9\textwidth]{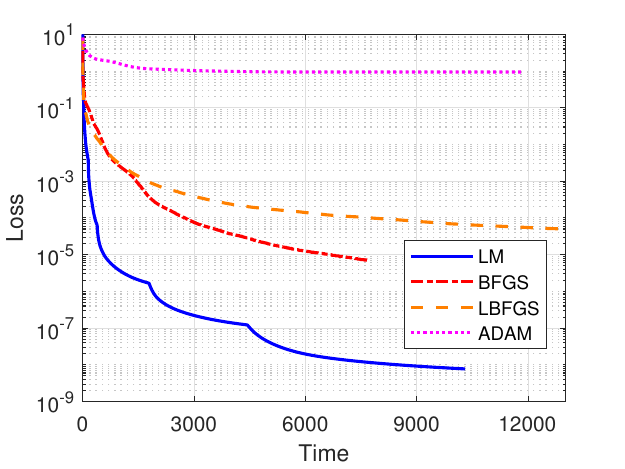}
		\caption{}
		\label{fig_forward_nls_loss_time}
	\end{subfigure}
	
	\caption{(a) Neural network solution of the nonlinear Schr\"odinger equation obtained using LM.
		(b) Solutions computed by LM and BFGS at selected time slices.
		(c) Absolute errors corresponding to the LM and BFGS solutions.
		(d) Training loss histories of LM, BFGS, L-BFGS, and Adam.}
	\label{fig_forward_nls}
\end{figure}

Figs.\ \ref{fig_forward_nls_solution1} and \ref{fig_forward_nls_solution2} show the spatiotemporal solution and comparisons at $t=0.59$, $0.79$, and $0.98$. The reference solution is computed using a fourth-order Runge--Kutta method with time step $\Delta t=(\pi/2)\times10^{-6}$ and a Fourier spectral discretization in space with $\Delta x=10/2^9$. The corresponding absolute errors are shown in Fig.\ \ref{fig_forward_nls_solution3}, where LM achieves errors approximately two orders of magnitude smaller than those of BFGS. Fig.\ \ref{fig_forward_nls_loss_time} shows the loss as a function of computational time, with LM consistently attaining lower loss values than the other methods.
To further quantify the solution accuracy, the relative $L_2$ errors are computed on a grid of $2^9$ spatial points and 1001 temporal points, yielding $2.1\times10^{-5}$ for LM and $2.7\times10^{-3}$ for BFGS. These results further confirm the higher accuracy of LM for this problem.

\begin{table}
	\small
	\centering
	\caption{Comparison of loss values, relative $L_2$ errors, and computational times for different network sizes obtained using LM for the nonlinear Schr\"odinger equation.}
	\label{tab_forward_nls_1}
	\begin{tabular}{|l|c|c|c|}
		\hline
		Network                  &       Loss        & Relative $L_2$ error & Time \\ \hline
		$2\times$NN$(3,5,5,1)$   & $8.5\cdot10^{-3}$ &  $9.1\cdot10^{-2}$   & 424  \\ \hline
		$2\times$NN$(3,10,10,1)$ & $6.6\cdot10^{-5}$ &  $3.4\cdot10^{-3}$   & 1198  \\ \hline
		$2\times$NN$(3,15,15,1)$ & $4.3\cdot10^{-7}$ &  $2.4\cdot10^{-4}$   & 3250  \\ \hline
		$2\times$NN$(3,20,20,1)$ & $2.4\cdot10^{-8}$ &  $4.3\cdot10^{-5}$   & 4992 \\ \hline
		$2\times$NN$(3,25,25,1)$ & $7.7\cdot10^{-9}$ &  $2.1\cdot10^{-5}$   & 8856 \\ \hline
	\end{tabular}
\end{table}

Tab.\ \ref{tab_forward_nls_1} summarizes the results for different network sizes trained using LM. Even the smaller $2\times$NN$(3,10,10,1)$ and $2\times$NN$(3,15,15,1)$ configurations achieve good accuracy, indicating that relatively small shallow architectures can be effective for this problem.

\subsection{Inverse problem: Allen--Cahn equation}

We next consider the application of PINNs to inverse problems using the Allen--Cahn equation \cite{raissi2019physics},
\begin{equation} \label{eq_allen_cahn1}
	\begin{split}
		& u_t - 0.0001 u_{xx} + 5 u^3 - 5 u = 0, \quad x \in [-1,1], \quad t \in [0,1], \\
		& u(x, 0) = x^2 \cos(\pi x), \\
		& u(-1, t) = u(1, t) = -1.
	\end{split}
\end{equation}
For the inverse problem, we consider the generalized form
\begin{equation} \label{eq_allen_cahn2}
	u_t + \lambda_1 u_{xx} + \lambda_2 u^3 + \lambda_3 u = 0
\end{equation}
where $\lambda_1$, $\lambda_2$, and $\lambda_3$ are unknown parameters to be identified simultaneously with the solution $u$.

The initial and boundary conditions are enforced through
\begin{equation} \label{eq_allen_cahn_enforce}
	u = \tilde{u} (1 - x^2) t  + x^2 \cos(\pi x).
\end{equation}
The reference solution is computed on $t\in[0,1]$ using a fourth-order Runge--Kutta method with time step $\Delta t=10^{-6}$, while the spatial derivatives are approximated using a central finite-difference scheme with $\Delta x=2\times10^{-4}$. The solution is then sampled on a uniform $1001\times1001$ grid with $\Delta x=2\times10^{-3}$ and $\Delta t=10^{-3}$. From this grid, $n=2000$ points are randomly selected without replacement and fixed throughout training, with the corresponding reference values used as observational data. The true parameters are $\lambda_1=-0.0001$, $\lambda_2=5$, and $\lambda_3=-5$.
Parameter identification is performed using NN$(2,25,25,1)$, giving 754 trainable parameters including the network weights, biases, and the unknown coefficients $\lambda_1$, $\lambda_2$, and $\lambda_3$. All three coefficients are initialized at zero.

\begin{figure}
	\centering
	
	\begin{subfigure}[b]{\textwidth}
		\centering
		\includegraphics[width=0.9\textwidth]{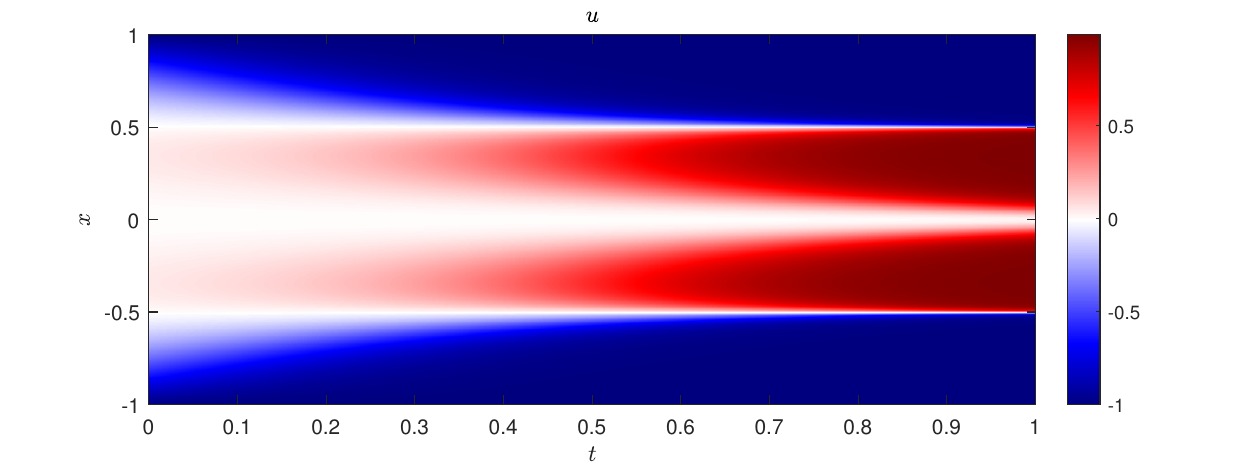}
		\caption{}
		\label{fig_inverse_allen_cahn_solution1}
	\end{subfigure}
	
	\begin{subfigure}[b]{\textwidth}
		\centering
		\includegraphics[width=\textwidth]{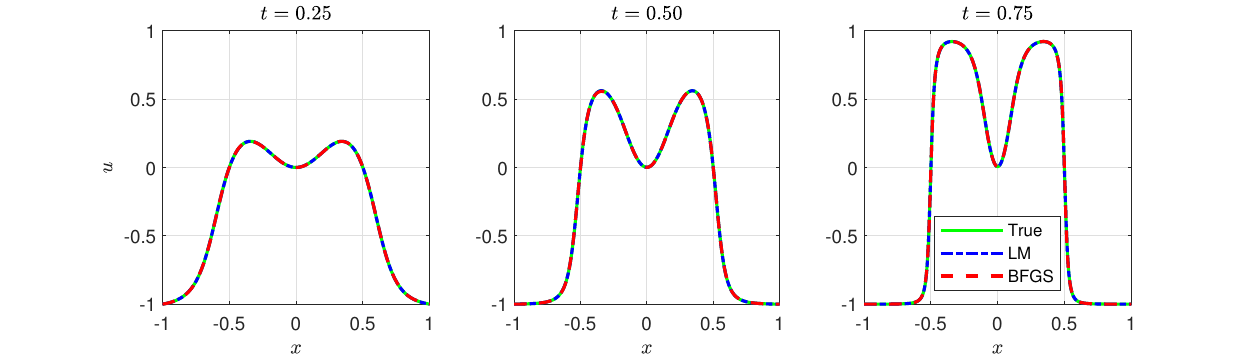}
		\caption{}
		\label{fig_inverse_allen_cahn_solution2}
	\end{subfigure}
	
	\begin{subfigure}[b]{\textwidth}
		\centering
		\includegraphics[width=\textwidth]{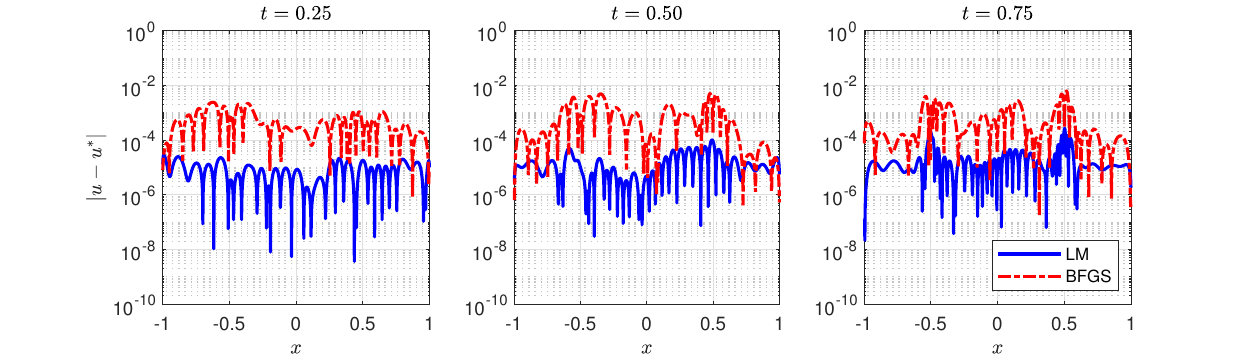}
		\caption{}
		\label{fig_inverse_allen_cahn_solution3}
	\end{subfigure}
	
	\begin{subfigure}[b]{0.49\textwidth}
		\centering
		\includegraphics[width=0.9\textwidth]{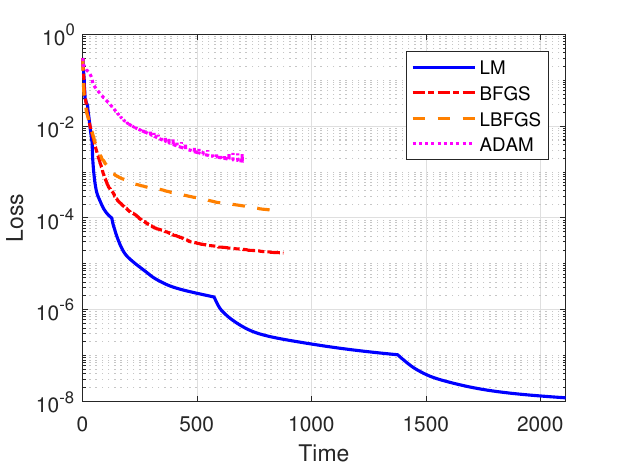}
		\caption{}
		\label{fig_inverse_allen_cahn_loss_time}
	\end{subfigure}
	\begin{subfigure}[b]{0.49\textwidth}
		\centering
		\includegraphics[width=0.9\textwidth]{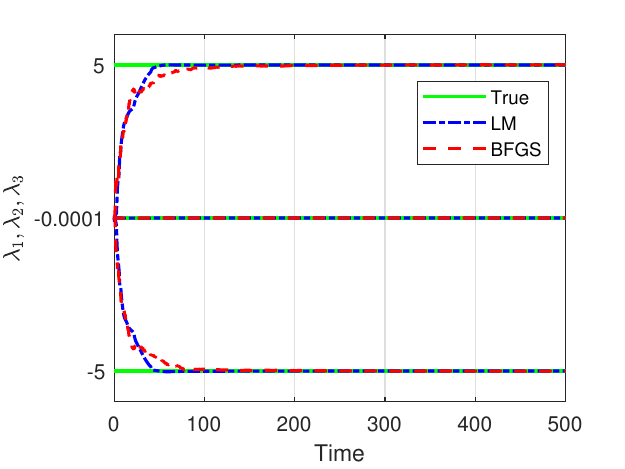}
		\caption{}
		\label{fig_inverse_allen_cahn_lambda}
	\end{subfigure}
	
	\caption{(a) Neural network solution of the Allen--Cahn equation obtained using LM.
		(b) Solutions computed by LM and BFGS at selected time slices.
		(c) Absolute errors corresponding to the LM and BFGS solutions.
		(d) Training loss histories of LM, BFGS, L-BFGS, and Adam.
		(e) Evolution of the estimated parameters during during the first 500 seconds.}
	\label{fig_inverse_allen_cahn}
\end{figure}

Fig.\ \ref{fig_inverse_allen_cahn_solution1} shows the full neural network solution obtained using LM. In Fig.\ \ref{fig_inverse_allen_cahn_solution2}, the LM and BFGS solutions at $t=0.25$, $0.50$, and $0.75$ are compared with the reference solution. The corresponding absolute errors are shown in Fig.\ \ref{fig_inverse_allen_cahn_solution3}. Fig.\ \ref{fig_inverse_allen_cahn_loss_time} shows the training loss as a function of computational time. Although LM requires more time to reach its final loss of $1.2\times10^{-8}$, it consistently attains lower loss values than the other methods over the considered time range.
Fig.\ \ref{fig_inverse_allen_cahn_lambda} shows the evolution of the estimated coefficients $\lambda_1$, $\lambda_2$, and $\lambda_3$ during the first 500 seconds, showing that LM converges to the true values faster than BFGS.

\begin{table}
	\small
	\centering
	\caption{Identified PDEs obtained using LM and BFGS for the Allen–Cahn equation.}
	\label{tab_inverse_allen_cahn_1}
	\begin{tabular}{|l|c|}
		\hline
		Method &                        PDE                        \\ \hline
		True   & $u_t - 0.000100000u_{xx} + 5.00000u^3 - 5.00000u$ \\ \hline
		LM     & $u_t - 0.000100003u_{xx} + 4.99983u^3 - 4.99994u$ \\ \hline
		BFGS   & $u_t - 0.000144231u_{xx} + 5.00372u^3 - 5.00402u$ \\ \hline
	\end{tabular}
\end{table}

\begin{table}
	\small
	\centering
	\caption{Comparison of loss values, absolute percentage errors (APEs), and computational times obtained using LM and BFGS for the Allen–Cahn equation.}
	\label{tab_inverse_allen_cahn_2}
	\begin{tabular}{|l|c|c|c|c|c|}
		\hline
		Method &       Loss        & APE $\lambda_1$ & APE $\lambda_2$ & APE $\lambda_3$ & Time \\ \hline
		LM     & $1.2\cdot10^{-8}$ &     0.0027      &     0.0034      &     0.0011      & 1677 \\ \hline
		BFGS   & $1.7\cdot10^{-5}$ &     44.2306     &     0.0744      &     0.0804      & 834  \\ \hline
	\end{tabular}
\end{table}

\begin{table}
	\small
	\centering
	\caption{Comparison of loss values, absolute percentage errors (APEs), and computational times for different network sizes obtained using LM for the Allen–Cahn equation.}
	\label{tab_inverse_allen_cahn_3}
	\begin{tabular}{|l|c|c|c|c|c|}
		\hline
		Network         &       Loss        & APE $\lambda_1$ & APE $\lambda_2$ & APE $\lambda_3$ & Time \\ \hline
		NN$(2,5,5,1)$   & $1.5\cdot10^{-2}$ &    1618.1627    &     25.4738     &     22.5804     &  40  \\ \hline
		NN$(2,10,10,1)$ & $5.7\cdot10^{-4}$ &     72.1185     &     2.6754      &     2.2376      & 110 \\ \hline
		NN$(2,15,15,1)$ & $8.8\cdot10^{-7}$ &     0.3691      &     0.0033      &     0.0036      & 244  \\ \hline
		NN$(2,20,20,1)$ & $8.1\cdot10^{-8}$ &     0.6517      &     0.0082      &     0.0016      & 684  \\ \hline
		NN$(2,25,25,1)$ & $1.2\cdot10^{-8}$ &     0.0027      &     0.0034      &     0.0011      & 1677 \\ \hline
	\end{tabular}
\end{table}

Tabs.\ \ref{tab_inverse_allen_cahn_1} and \ref{tab_inverse_allen_cahn_2} summarize the identified equations and the corresponding absolute percentage errors (APEs). LM identifies all three coefficients more accurately than BFGS, demonstrating the effectiveness of shallow PINNs trained with LM for parameter identification in nonlinear PDEs.
Tab.\ \ref{tab_inverse_allen_cahn_3} compares the performance of different network sizes trained using LM. Increasing the number of neurons generally improves the loss and parameter estimates. Notably, NN$(2,15,15,1)$ already achieves APEs below $1\%$ for all coefficients and outperforms BFGS in both accuracy and computational time, indicating that a moderately sized shallow network can provide an effective balance between efficiency and accuracy.

\begin{table}[h]
	\small
	\centering
	\caption{Identified PDEs obtained using LM for different noise settings.}
	\label{tab_inverse_allen_cahn_4}
	\begin{tabular}{|c|c|}
		\hline
		Noise &                        PDE                        \\ \hline
		0\%   & $u_t - 0.000100003u_{xx} + 4.99983u^3 - 4.99994u$ \\ \hline
		1\%   &  $u_t -0.000094419u_{xx} + 4.99628u^3 -4.99628u$  \\ \hline
		2\%   &  $u_t -0.000101316u_{xx} + 4.99665u^3 -4.99676u$  \\ \hline
		3\%   &  $u_t -0.000102621u_{xx} + 5.00054u^3 -4.99968u$  \\ \hline
		4\%   &  $u_t -0.000091892u_{xx} + 4.99998u^3 -4.99901u$  \\ \hline
		5\%   &  $u_t -0.000093913u_{xx} + 5.01459u^3 -5.01438u$  \\ \hline
		6\%   &  $u_t -0.000116742u_{xx} + 4.99733u^3 -5.00159u$  \\ \hline
		7\%   &  $u_t -0.000125068u_{xx} + 5.00743u^3 -5.01171u$  \\ \hline
		8\%   &  $u_t -0.000052940u_{xx} + 4.99270u^3 -4.98783u$  \\ \hline
		9\%   &  $u_t -0.000267716u_{xx} + 4.97542u^3 -4.99050u$  \\ \hline
		10\%  &  $u_t -0.000096346u_{xx} + 4.98809u^3 -4.98886u$  \\ \hline
	\end{tabular}
\end{table}

\begin{table}[h]
	\small
	\centering
	\caption{Identified PDEs obtained using LM for different initial guess $(\lambda_1, \lambda_2, \lambda_3)$.}
	\label{tab_inverse_allen_cahn_5}
	\begin{tabular}{|r|r|r|c|}
		\hline
		$\lambda_1$ & $\lambda_2$ & $\lambda_3$ &                        PDE                        \\ \hline
		  $0.00000$ &   $0.00000$ &   $0.00000$ & $u_t - 0.000100003u_{xx} + 4.99983u^3 - 4.99994u$ \\ \hline
		  $0.00054$ &   $5.36834$ &   $3.74416$ &  $u_t -0.000099667u_{xx} + 4.99973u^3 -4.99991u$  \\ \hline
		 $-0.00018$ &   $3.90923$ &  $-6.53991$ &  $u_t -0.000099913u_{xx} + 4.99987u^3 -4.99998u$  \\ \hline
		  $0.00069$ &  $-4.57194$ &   $1.58149$ &  $u_t -0.000099566u_{xx} + 4.99987u^3 -4.99993u$  \\ \hline
		 $-0.00059$ &  $-0.22215$ &  $-8.03462$ &  $u_t -0.000099788u_{xx} + 4.99984u^3 -4.99998u$  \\ \hline
	\end{tabular}
\end{table}

To assess robustness to observational noise, uniformly distributed multiplicative noise is added to the reference solution according to
\begin{equation*}
	\tilde{u}_i=u_i(1+\alpha\xi_i),\qquad \xi_i\sim\mathcal{U}(-1,1),
\end{equation*}
where $\alpha$ denotes the noise level. Thus, $\alpha=0.01$ corresponds to a relative perturbation within $\pm1\%$ of each reference value. Noise levels from $1\%$ to $10\%$ are considered. As shown in Tab.\ \ref{tab_inverse_allen_cahn_4}, LM recovers the governing equation accurately even at relatively high noise levels. The identified coefficients remain close to their true values, although the diffusion coefficient $\lambda_1$ shows greater sensitivity to noise.

The sensitivity to initialization is examined in Tab.\ \ref{tab_inverse_allen_cahn_5}. Despite substantially different initial guesses, the recovered coefficients are highly consistent and remain close to the true values. Together, these experiments provide numerical evidence that the proposed shallow LM-PINN is robust to both observational noise and parameter initialization for the inverse Allen--Cahn problem.

\subsection{Inverse problem: three-dimensional Bratu equation}

As the final example, we consider a time-independent inverse problem governed by the three-dimensional Bratu equation. Unlike the previous examples, which involve a single spatial variable, this problem illustrates the applicability of shallow PINNs to higher-dimensional spatial domains. The three-dimensional Bratu equation is given by \cite{shahab2025finite}
\begin{equation} \label{eq_bratu1}
	\begin{split}
		& \Delta u + C \exp(u) = 0, \quad (x, y, z) \in \Omega \\
		& u(x, y, z) = 0, \quad (x, y, z) \in \partial\Omega
	\end{split}
\end{equation}
where $\Omega = [0,1]^3$ and $\Delta u = u_{xx} + u_{yy} + u_{zz}$.
For the inverse problem, we consider
\begin{equation}
	\Delta u + \lambda_1 \exp(\lambda_2 u) = 0
\end{equation}
where $\lambda_1$ and $\lambda_2$ are unknown parameters to be identified together with the solution $u$.

The boundary conditions are enforced through
\begin{equation} \label{eq_bratu_enforce}
	u = \tilde{u} (x-x^2) (y-y^2) (z-z^2).
\end{equation}
For training, $n=10000$ collocation points are randomly selected without replacement from a uniform grid over the computational domain. The reference solution is computed using the symmetric finite difference method (SFDM) \cite{shahab2025finite} with 300 uniform subintervals, corresponding to 301 grid points in each spatial direction. The selected collocation points are fixed throughout training, and the corresponding reference values are used as observational data. The true parameters are $\lambda_1=C=2$ and $\lambda_2=1$.
Parameter identification is performed using NN$(3,25,25,1)$, resulting in 778 trainable parameters, including the network weights, biases, and the unknown coefficients $\lambda_1$ and $\lambda_2$. Both coefficients are initialized at zero.

\begin{figure}
	\centering
	
	\begin{subfigure}[b]{0.49\textwidth}
		\centering
		\includegraphics[width=0.9\textwidth]{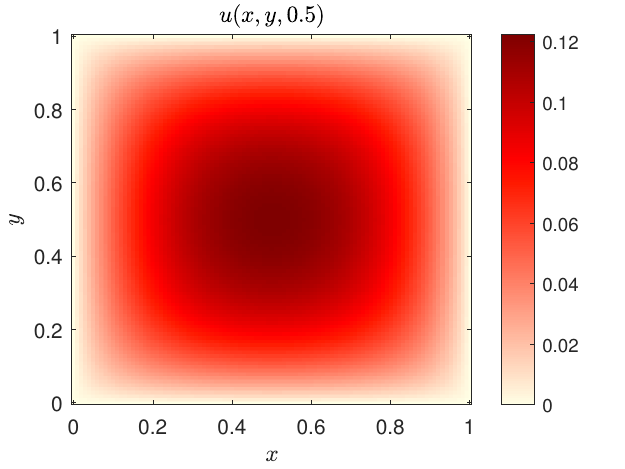}
		\caption{}
		\label{fig_inverse_bratu_solution1}
	\end{subfigure}
	
	\begin{subfigure}[b]{\textwidth}
		\centering
		\includegraphics[width=\textwidth]{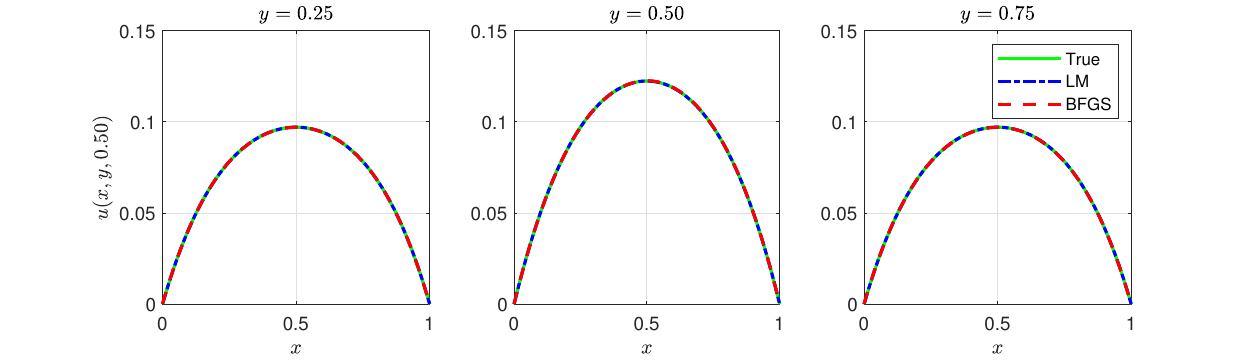}
		\caption{}
		\label{fig_inverse_bratu_solution2}
	\end{subfigure}
	
	\begin{subfigure}[b]{\textwidth}
		\centering
		\includegraphics[width=\textwidth]{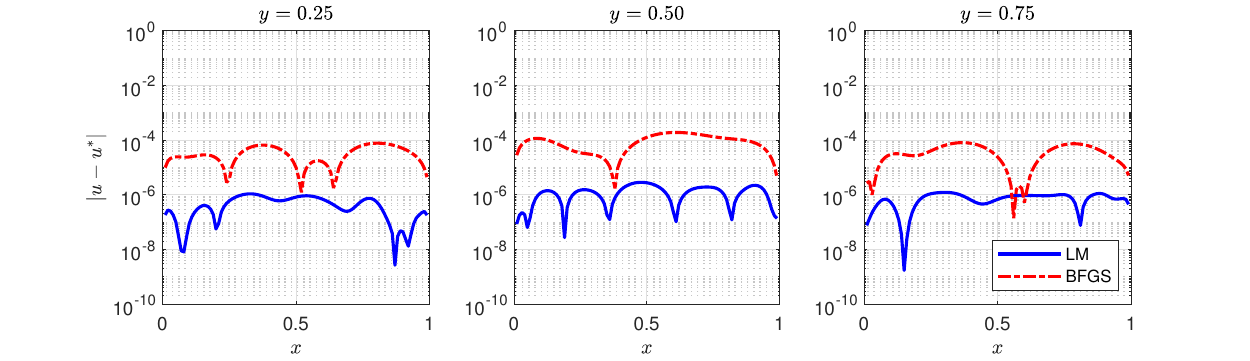}
		\caption{}
		\label{fig_inverse_bratu_solution3}
	\end{subfigure}
	
	\begin{subfigure}[b]{0.49\textwidth}
		\centering
		\includegraphics[width=0.9\textwidth]{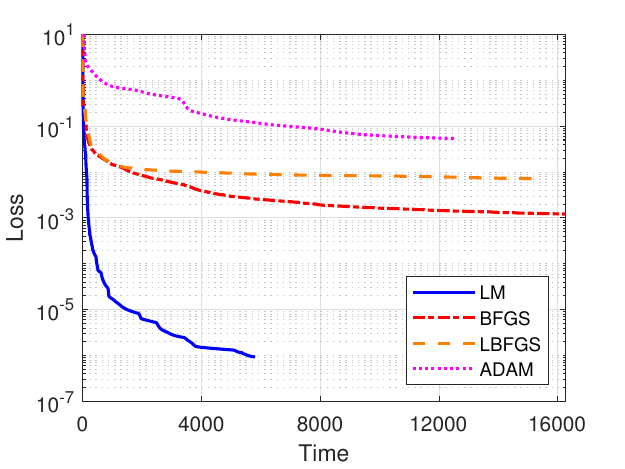}
		\caption{}
		\label{fig_inverse_bratu_loss_time}
	\end{subfigure}
	\begin{subfigure}[b]{0.49\textwidth}
		\centering
		\includegraphics[width=0.9\textwidth]{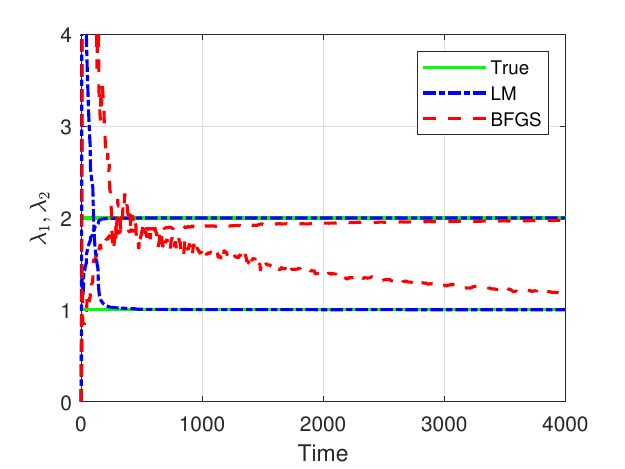}
		\caption{}
		\label{fig_inverse_bratu_lambda}
	\end{subfigure}
	
	\caption{(a) Neural network solution of the three-dimensional Bratu equation obtained using LM.
		(b) Solutions computed by LM and BFGS at selected time slices.
		(c) Absolute errors corresponding to the LM and BFGS solutions.
		(d) Training loss histories of LM, BFGS, L-BFGS, and Adam.
		(e) Evolution of the estimated parameters during the first 4000 seconds.}
	\label{fig_inverse_bratu}
\end{figure}

The neural network solution obtained using LM is shown in Fig.\ \ref{fig_inverse_bratu_solution1} on the plane $z=0.5$. In Fig.\ \ref{fig_inverse_bratu_solution2}, the LM and BFGS solutions at $z=0.5$ and $y=0.25$, $0.50$, and $0.75$ are compared with the SFDM reference solution. The corresponding absolute errors are shown in Fig.\ \ref{fig_inverse_bratu_solution3}. Fig.\ \ref{fig_inverse_bratu_loss_time} shows the training loss as a function of computational time, with LM consistently attaining lower loss values than the other methods over the considered time range. Fig.\ \ref{fig_inverse_bratu_lambda} shows the evolution of the estimated coefficients $\lambda_1$ and $\lambda_2$ during the first 4000 seconds, demonstrating that LM converges to the true values more rapidly and accurately than BFGS.

\begin{table}[t]
	\small
	\centering
	\caption{Identified PDEs obtained using LM and BFGS for the three-dimensional Bratu equation.}
	\label{tab_inverse_bratu_1}
	\begin{tabular}{|l|c|}
		\hline
		Method &                 PDE                  \\ \hline
		True   & $\Delta u + 2.00000 \exp(1.00000 u)$ \\ \hline
		LM     & $\Delta u + 1.99992 \exp(1.00040 u)$ \\ \hline
		BFGS   & $\Delta u + 1.99144 \exp(1.05985 u)$ \\ \hline
	\end{tabular}
\end{table}

\begin{table}[t]
	\small
	\centering
	\caption{Comparison of loss values, absolute percentage errors (APEs), and computational times obtained using LM and BFGS for the three-dimensional Bratu equation.}
	\label{tab_inverse_bratu_2}
	\begin{tabular}{|l|c|c|c|c|}
		\hline
		Method &       Loss        & APE $\lambda_1$ & APE $\lambda_2$ & Time  \\ \hline
		LM     & $9.4\cdot10^{-7}$ &     0.0041      &     0.0401      & 5074  \\ \hline
		BFGS   & $1.2\cdot10^{-3}$ &     0.4281      &     5.9845      & 14465 \\ \hline
	\end{tabular}
\end{table}

\begin{table}[t]
	\small
	\centering
	\caption{Comparison of loss values, absolute percentage errors (APEs), and computational times for different network sizes obtained using LM for the three-dimensional Bratu equation.}
	\label{tab_inverse_bratu_3}
	\begin{tabular}{|l|c|c|c|c|}
		\hline
		Network         &       Loss        & APE $\lambda_1$ & APE $\lambda_2$ & Time \\ \hline
		NN$(3,5,5,1)$   & $2.7\cdot10^{-2}$ &     5.1426      &     71.8921     & 274  \\ \hline
		NN$(3,10,10,1)$ & $2.9\cdot10^{-3}$ &     0.9318      &     12.8455     & 725  \\ \hline
		NN$(3,15,15,1)$ & $2.0\cdot10^{-4}$ &     0.0560      &     0.8668      & 2161 \\ \hline
		NN$(3,20,20,1)$ & $4.9\cdot10^{-5}$ &     0.0212      &     0.3051      & 3627 \\ \hline
		NN$(3,25,25,1)$ & $9.4\cdot10^{-7}$ &     0.0041      &     0.0401      & 5074 \\ \hline
	\end{tabular}
\end{table}

Tabs.\ \ref{tab_inverse_bratu_1} and \ref{tab_inverse_bratu_2} summarize the identified equations and corresponding APEs. LM identifies both coefficients substantially more accurately than BFGS while also requiring less computational time, further demonstrating the effectiveness of shallow PINNs trained with LM for nonlinear PDE parameter identification.
Tab.\ \ref{tab_inverse_bratu_3} compares different network sizes trained using LM. Increasing the number of neurons generally improves both the loss and parameter estimates. Notably, NN$(3,15,15,1)$ already achieves APEs below $1\%$ for both coefficients and outperforms BFGS in both accuracy and computational time, indicating that a moderately sized shallow network can provide an effective balance between efficiency and accuracy.

\begin{table}[t]
	\small
	\centering
	\caption{Identified PDEs obtained using LM for different noise settings.}
	\label{tab_inverse_bratu_4}
	\begin{tabular}{|c|c|}
		\hline
		Noise &                 PDE                  \\ \hline
		 0\%  & $\Delta u + 1.99992 \exp(1.00040 u)$ \\ \hline
		 1\%  & $\Delta u + 1.99733 \exp(1.01903 u)$ \\ \hline
		 2\%  & $\Delta u + 2.00184 \exp(0.98412 u)$ \\ \hline
		 3\%  & $\Delta u + 1.99398 \exp(1.04735 u)$ \\ \hline
		 4\%  & $\Delta u + 1.98806 \exp(1.08270 u)$ \\ \hline
		 5\%  & $\Delta u + 1.98904 \exp(1.07992 u)$ \\ \hline
		 6\%  & $\Delta u + 2.00316 \exp(0.97100 u)$ \\ \hline
		 7\%  & $\Delta u + 2.00725 \exp(0.94867 u)$ \\ \hline
		 8\%  & $\Delta u + 2.00997 \exp(0.91204 u)$ \\ \hline
		 9\%  & $\Delta u + 1.99150 \exp(1.06681 u)$ \\ \hline
		10\%  & $\Delta u + 2.00433 \exp(0.97587 u)$ \\ \hline
	\end{tabular}
\end{table}

\begin{table}[h]
	\small
	\centering
	\caption{Identified PDEs obtained using LM for different initial guess.}
	\label{tab_inverse_bratu_5}
	\begin{tabular}{|c|c|}
		\hline
		Initial Guess $(\lambda_1, \lambda_2)$ &                 PDE                  \\ \hline
		               $(0, 0)$                & $\Delta u + 1.99992 \exp(1.00040 u)$ \\ \hline
		         $(6.51790, 3.82647)$          & $\Delta u + 1.99952 \exp(1.00310 u)$ \\ \hline
		         $(7.80156, 4.88854)$          & $\Delta u + 1.99991 \exp(1.00039 u)$ \\ \hline
		         $(3.99433, 5.20323)$          & $\Delta u + 1.99949 \exp(1.00318 u)$ \\ \hline
		         $(2.16443, 7.96856)$          & $\Delta u + 1.99979 \exp(1.00121 u)$ \\ \hline
	\end{tabular}
\end{table}

To assess robustness to observational noise, we use the same noise-generation procedure as in the inverse Allen--Cahn problem. As shown in Tab.\ \ref{tab_inverse_bratu_4}, the identified coefficients remain reasonably close to their reference values, $\lambda_1=2$ and $\lambda_2=1$, even at a noise level of $10\%$. Although $\lambda_2$ is more sensitive to increasing noise than $\lambda_1$, the governing PDE is consistently recovered. Tab.\ \ref{tab_inverse_bratu_5} examines sensitivity to initialization. Despite substantially different initial guesses, LM consistently converges to coefficients close to the reference values. These results provide numerical evidence that parameter identification for the three-dimensional Bratu problem is robust to observational noise and relatively insensitive to initialization.

\subsection{Discussion on the Damping Parameter and Computational Cost}

This section examines the role of the LM damping parameter $\mu$ and the computational cost of the proposed approach. Unless otherwise stated, all experiments are conducted for the Burgers equation in Eq.\ \eqref{eq_burgers1} using the network architecture $\mathrm{NN}(2,25,25,1)$. An additional experiment with deeper networks is presented at the end of the section to further investigate the effect of network depth and provide a comparison with L-BFGS.

The Gauss--Newton update is given by
\begin{equation}
	\label{eq_Gauss_Newton_update}
	\Delta W = (J^T J)^{-1}J^TF.
\end{equation}
The main difference between Gauss--Newton and LM is the damping term $\mu I$ added to the normal equations, as shown in Eq.\ \eqref{eq_LM_update}. At the first iteration of the Burgers problem, the Gauss--Newton formulation yields a 2-norm condition number of
$\kappa(J^TJ)=2.16\times10^{22}$, indicating an extremely ill-conditioned linear system. Adding the damping term with $\mu=100$ reduces the condition number to
$\kappa(J^TJ+\mu I)=1.13\times10^{3}$.
Thus, damping substantially improves the conditioning and numerical stability of the linear system, illustrating an important practical advantage of LM over the undamped Gauss--Newton method for the ill-conditioned nonlinear least-squares problems encountered in PINN training.

\begin{figure}
	\centering
	
	\begin{subfigure}[b]{0.49\textwidth}
		\centering
		\includegraphics[width=0.9\textwidth]{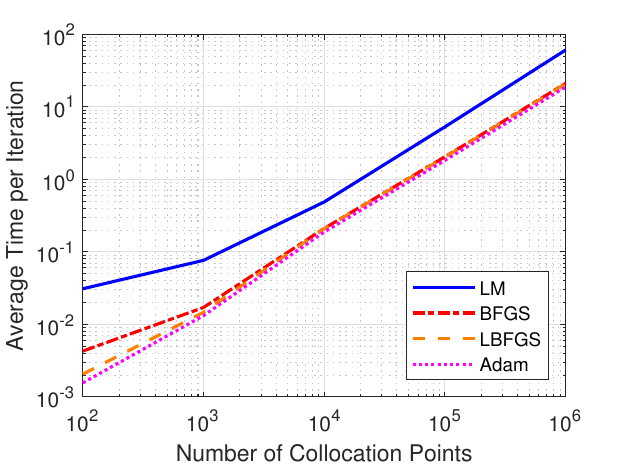}
		\caption{}
		\label{fig_forward_burgers_average_time}
	\end{subfigure}
	\begin{subfigure}[b]{0.49\textwidth}
		\centering
		\includegraphics[width=0.9\textwidth]{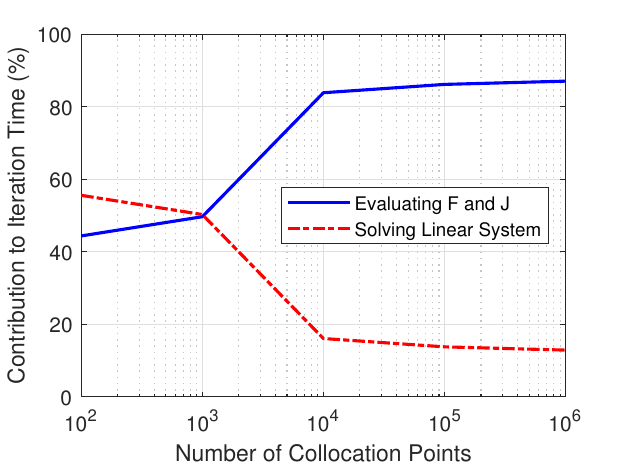}
		\caption{}
		\label{fig_forward_burgers_contribution_time}
	\end{subfigure}
	
	\caption{(a) Average time per iteration for different numbers of collocation points. (b) Percentage contribution of evaluating $F$ and $J$ and solving the LM normal equations to the total iteration time.}
	\label{fig_forward_burgers_time}
\end{figure}

We next examine the computational cost of LM. Fig.\ \ref{fig_forward_burgers_average_time} compares the average iteration times of LM, BFGS, L-BFGS, and Adam for $10^2$, $10^3$, $10^4$, $10^5$, and $10^6$ collocation points. Although each LM iteration is more expensive than an iteration of the other optimizers, the computational cost of all methods grows approximately linearly with the number of collocation points. For $10^4$ collocation points and above, the rate of increase for LM is comparable to those of BFGS, L-BFGS, and Adam.

Fig.\ \ref{fig_forward_burgers_contribution_time} further decomposes the LM iteration time into residual and Jacobian evaluation and solution of the LM normal equations. For small numbers of collocation points, the linear solve accounts for more than $50\%$ of the total iteration time. Its relative contribution decreases substantially as the number of collocation points increases, falling below $20\%$ for $10^4$ or more points and to approximately $13\%$ for $10^6$ points. Residual and Jacobian evaluation therefore becomes the dominant computational cost for large problems. This indicates that, for the PINNs considered in this work, the computational cost of LM is dominated by residual and Jacobian evaluation rather than by solving the LM linear system.

\begin{figure}
	\centering
	
	\includegraphics[width=0.49\textwidth]{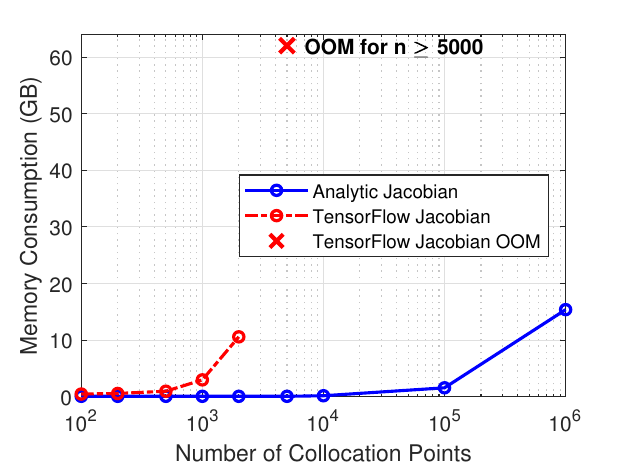}	
	\caption{Memory consumption for computing the residual Jacobian using the proposed analytical formulation and TensorFlow for different numbers of collocation points. TensorFlow runs out of memory for $n\geq5000$.}
	\label{fig_forward_burgers_memory}
\end{figure}

The memory requirements of the residual Jacobian computation are examined in Fig.\ \ref{fig_forward_burgers_memory}. The proposed analytical formulation exhibits substantially slower memory growth than the automatic-differentiation-based implementation in TensorFlow. While the analytical formulation remains feasible for up to $10^6$ collocation points, the TensorFlow implementation runs out of memory at approximately $5000$ points under the tested configuration. These results highlight the practical advantage of explicitly deriving the residual Jacobian when applying LM to PINNs with large numbers of collocation points.

\begin{table}
	\small
	\centering
	\caption{Loss values obtained for the Burgers equation using L-BFGS with different numbers of hidden layers.}
	\label{tab_forward_burgers_2}
	\begin{tabular}{|l|c|c|c|c|}
		\hline
		     &     3 Layers      &     4 Layers      &     5 Layers      &     6 Layers      \\ \hline
		Loss & $1.9\cdot10^{-4}$ & $9.6\cdot10^{-6}$ & $7.2\cdot10^{-6}$ & $3.6\cdot10^{-6}$ \\ \hline
		\hline
		     &     7 Layers      &     8 Layers      &     9 Layers      &     10 Layers     \\ \hline
		Loss & $4.3\cdot10^{-6}$ & $1.2\cdot10^{-5}$ & $2.2\cdot10^{-6}$ & $3.7\cdot10^{-6}$ \\ \hline
	\end{tabular}
\end{table}

Finally, we examine whether increasing network depth is necessary for achieving accurate PINN solutions. Tab.\ \ref{tab_forward_burgers_2} reports the losses obtained using L-BFGS for networks with three to ten hidden layers, each containing 25 neurons. These experiments are implemented in PyTorch. For comparison, the shallow two-hidden-layer network used in this work achieves a loss of $8.1\times10^{-8}$ with LM, which is lower than all losses obtained by the deeper L-BFGS networks. Although this is not a direct optimizer comparison under identical architectures, it provides additional numerical evidence that greater network depth is not necessarily required for high accuracy when an effective second-order optimization method is employed.

\section{Conclusion}

In this work, we investigated shallow PINNs for solving nonlinear PDEs by reformulating PINN training as a nonlinear least-squares problem and solving it using the Levenberg--Marquardt (LM) algorithm. The results demonstrate that shallow PINNs, when combined with an effective second-order optimization method, can achieve high accuracy without relying on deep network architectures.

We derived exact analytical expressions for neural-network derivatives with respect to the input variables. These expressions play a central role in constructing the proposed PINN framework (see Fig.\ \ref{fig_pinns}), as they explicitly reveal the relationships between the network output and its spatial and temporal derivatives while providing a transparent representation of the parameter-sharing structure of PINNs. We also derived explicit formulations for the Jacobian required by LM and established its relationship to the gradient of the PINN loss used in gradient-based and quasi-Newton methods.

Numerical experiments were conducted for forward and inverse problems involving the Burgers, Schr\"odinger, Allen--Cahn, and three-dimensional Bratu equations. These problems cover time-dependent and time-independent PDEs, different nonlinearities, multiple spatial dimensions, and both real- and complex-valued solution representations. Using the same shallow architecture $\mathrm{NN}(a,25,25,1)$, LM consistently achieved substantially lower loss values than BFGS, L-BFGS, and Adam at comparable computational times. For the forward problems, LM produced highly accurate solution approximations, with relative $L_2$ errors substantially smaller than those obtained by the other methods. For the inverse problems, LM accurately recovered the unknown PDE parameters, demonstrating its effectiveness for both solution approximation and parameter identification.

Comparison based on computational time is particularly important because each LM iteration is generally more expensive than an iteration of a gradient-based or quasi-Newton method. Nevertheless, the faster convergence of LM compensates for this additional per-iteration cost, allowing it to attain lower loss values within comparable computational times. Additional experiments with deeper networks trained using L-BFGS further support this observation. For the Burgers equation, increasing the network depth from three to ten hidden layers did not yield a lower loss than the shallow two-hidden-layer network trained with LM. Thus, in the tested setting, the shallow LM-PINN achieved higher accuracy with substantially fewer parameters, suggesting that network architecture and optimization strategy should be considered jointly.

The analysis of the LM damping parameter provides further insight into the numerical behavior of the proposed method. For the Burgers equation, the Gauss--Newton system was extremely ill-conditioned, whereas the LM damping term substantially reduced the condition number and improved the numerical stability of the resulting linear system. This observation helps explain the robustness of LM in the PINN optimization problems considered here and highlights the importance of damping in distinguishing LM from the undamped Gauss--Newton method.

The computational and memory analyses further clarify the scalability of the proposed approach. As the number of collocation points increases, the dominant computational cost of LM shifts toward residual and Jacobian evaluation, while the relative cost of solving the LM linear system decreases. Thus, for large PINN problems, reducing the cost of residual and Jacobian evaluation may be more beneficial than optimizing the linear solver alone. Moreover, the explicit analytical Jacobian provides a substantial memory advantage over automatic differentiation, indicating that analytical differentiation can be particularly advantageous for large-scale PINNs.

Overall, the results suggest that high accuracy in PINNs does not necessarily require deep architectures when effective optimization and accurate derivative information are available. We do not claim that shallow networks are universally superior to deep architectures. Instead, the appropriate network depth should be determined jointly by the characteristics of the mathematical problem, the desired solution representation, and the optimization method. In practice, shallow architectures such as $\mathrm{NN}(a,15,15,1)$ may provide a useful starting point, while additional neurons or hidden layers can be introduced when greater representational capacity is required, particularly for high-frequency or multiscale solutions \cite{wang2021eigenvector}.

Beyond its numerical performance, the present study provides an analytical perspective on the internal structure of PINNs. By explicitly deriving network derivatives, residuals, Jacobians, and their relationships to loss gradients, the proposed formulation offers a transparent alternative to treating automatic differentiation as a black-box operation. These analytical expressions are not specific to LM and may also support the implementation and analysis of other optimization strategies for PINNs.

Several directions remain for future research. Although the present study includes a three-dimensional PDE, further investigation is needed to determine whether the observed advantages of shallow LM-PINNs persist in substantially higher spatial dimensions, where computational and memory requirements can increase rapidly and more expressive architectures or alternative sampling strategies may be required. Stochastic Levenberg--Marquardt methods \cite{shahab2025physics} therefore represent a promising direction for reducing the computational cost of large-scale problems. In addition, formulations that reduce or avoid explicit Jacobian storage, such as the improved LM approaches proposed in \cite{wilamowski2010improved}, may help alleviate the memory requirements of large-scale implementations.

\section*{Acknowledgement}
This research is funded by Institut Teknologi Sepuluh Nopember (ITS) and managed under the Strategic Research Grant (SRG) Type C Scheme (Contract No. 1453/PKS/ITS/2026).
HS acknowledges support by Khalifa University through the Research \& Innovation Grants under project IDs KU-INT-RIG-2023-8474000617 and KU-INT-RIG-2024-8474000789.

\section*{Declaration of competing interest} 
The authors declare that they have no known competing financial interests or personal relationships that could have appeared to influence the work reported in this paper.

\section*{Declaration of generative AI and AI-assisted technologies in the writing process}
During the preparation of this work, the authors used ChatGPT to improve the language and readability. 
After using this tool, the authors reviewed and edited the content as needed and take full responsibility for the content of the publication.

\bibliographystyle{elsarticle-num} 
\bibliography{main_references}

\appendix

\section{Derivatives with respect to weights and biases} \label{sec_derivatives_weight}

In this section, we provide the formulas for computing $\frac{\partial \tilde{u}}{\partial w_r}$, $\frac{\partial \tilde{u}_t}{\partial w_r}$, $\frac{\partial \tilde{u}x}{\partial w_r}$, and $\frac{\partial \tilde{u}{xx}}{\partial w_r}$. To simplify and accelerate the computations, we adopt the vectorized form $X = [\textbf{x} \; \textbf{t}]$, as described in Section \ref{sec_derivatives_input}. We begin with the weights and biases closest to the output layer. The derivatives with respect to the bias $B_3$ are given by
\begin{equation}
	\begin{split}
		\frac{\partial \tilde{u}}{\partial B_3} & = \textbf{1} \\
		\frac{\partial \tilde{u}_t}{\partial B_3} & = \textbf{0} \\
		\frac{\partial \tilde{u}_x}{\partial B_3} & = \textbf{0} \\
		\frac{\partial \tilde{u}_{xx}}{\partial B_3} & = \textbf{0} .
	\end{split}
\end{equation}

The derivatives with respect to the weights $w_j \in W_3$, for $j = 1, \dots, m_2$, are given by
\begin{equation}
	\begin{split}
		\frac{\partial \tilde{u}}{\partial w_j} & = H_2(:,j) \\
		\frac{\partial \tilde{u}_t}{\partial w_j} & = H_{2,t}(:,j) \\
		\frac{\partial \tilde{u}_x}{\partial w_j} & = H_{2,x}(:,j) \\
		\frac{\partial \tilde{u}_{xx}}{\partial w_j} & = H_{2,xx}(:,j) .
	\end{split}
\end{equation}
In this case, we use the notation $(:, j)$ to represent all rows in the $j$-th column of a matrix. This is a common convention in MATLAB. We adopt this notation because many of the matrices already include indices.

The derivatives with respect to the biases $b_k \in B_2$, for $k=1,\dots,m_2$, are given by
\begin{equation}
	\begin{split}
		\frac{\partial \tilde{u}}{\partial b_k} & = \sigma'(h_2(:,k)) \cdot W_3(k) \\
		\frac{\partial \tilde{u}_t}{\partial b_k} & = \sigma''(h_2(:,k)) \cdot h_{2,t}(:,k) \cdot W_3(k) \\
		\frac{\partial \tilde{u}_x}{\partial b_k} & = \sigma''(h_2(:,k)) \cdot h_{2,x}(:,k) \cdot W_3(k) \\
		\frac{\partial \tilde{u}_{xx}}{\partial b_k} 
		& = \sigma'''(h_2(:,k)) \cdot h_{2,x}^2(:,k) \cdot W_3(k) \\
		& \;\;\;\; + \sigma''(h_2(:,k)) \cdot h_{2,xx}(:,k) \cdot W_3(k) .
	\end{split}
\end{equation}
Similar to Section \ref{sec_derivatives_input}, the dot operator $\cdot$ here denotes element-wise multiplication.

The derivatives with respect to the weights $w_{j,k} \in W_2$, for $j = 1, \dots, m_1$ and $k=1,\dots,m_2$, are given by
\begin{equation}
	\begin{split}
		\frac{\partial \tilde{u}}{\partial w_{j,k}} 
		& = H_1(:,j) \cdot \sigma'(h_2(:,k)) \cdot W_3(k) \\
		\frac{\partial \tilde{u}_t}{\partial w_{j,k}} 
		& = H_1(:,j) \cdot \sigma''(h_2(:,k)) \cdot h_{2,t}(:,k) \cdot W_3(k) \\
		& \;\;\;\; + H_{1,t}(:,j) \cdot \sigma'(h_2(:,k)) \cdot W_3(k) \\
		\frac{\partial \tilde{u}_x}{\partial w_{j,k}} 
		& = H_1(:,j) \cdot \sigma''(h_2(:,k)) \cdot h_{2,x}(:,k) \cdot W_3(k) \\
		& \;\;\;\; + H_{1,x}(:,j) \cdot \sigma'(h_2(:,k)) \cdot W_3(k) \\
		\frac{\partial \tilde{u}_{xx}}{\partial w_{j,k}} 
		& = H_1(:,j) \cdot \sigma'''(h_2(:,k)) \cdot h_{2,x}^2(:,k) \cdot W_3(k) \\
		& \;\;\;\; + 2 H_{1,x}(:,j) \cdot \sigma''(h_2(:,k)) \cdot h_{2,x}(:,k) \cdot W_3(k) \\
		& \;\;\;\; + H_1(:,j) \cdot \sigma''(h_2(:,k)) \cdot h_{2,xx}(:,k) \cdot W_3(k) \\
		& \;\;\;\; + H_{1,xx}(:,j) \cdot \sigma'(h_2(:,k)) \cdot W_3(k) .
	\end{split}
\end{equation}

The derivatives with respect to the biases $b_k \in B_1$, for $k=1,\dots,m_1$, are given by
\begin{equation}
	\begin{split}
		\frac{\partial \tilde{u}}{\partial b_k} 
		& = \sigma'(h_1(:,k)) \cdot ((\sigma'(h_2) \cdot W_3^T) W_2(k,:)^T) \\
		\frac{\partial \tilde{u}_t}{\partial b_k} 
		& = \sigma'(h_1(:,k)) \cdot ((\sigma''(h_2) \cdot h_{2,t} \cdot W_3^T) W_2(k,:)^T) \\
		& \;\;\;\; + \sigma''(h_1(:,k)) \cdot h_{1,t}(:,k) \cdot ((\sigma'(h_2) \cdot W_3^T) W_2(k,:)^T) \\
		\frac{\partial \tilde{u}_x}{\partial b_k} 
		& = \sigma'(h_1(:,k)) \cdot ((\sigma''(h_2) \cdot h_{2,x} \cdot W_3^T) W_2(k,:)^T) \\
		& \;\;\;\; + \sigma''(h_1(:,k)) \cdot h_{1,x}(:,k) \cdot ((\sigma'(h_2) \cdot W_3^T) W_2(k,:)^T) \\
		\frac{\partial \tilde{u}_{xx}}{\partial b_k} 
		& = \sigma'(h_1(:,k)) \cdot ((\sigma'''(h_2) \cdot h_{2,x}^2 \cdot W_3^T) W_2(k,:)^T) \\
		& \;\;\;\; + \sigma''(h_1(:,k)) \cdot h_{1,x}(:,k) \cdot ((2 \sigma''(h_2) \cdot h_{2,x} \cdot W_3^T) W_2(k,:)^T) \\
		& \;\;\;\; + \sigma'(h_1(:,k)) \cdot ((\sigma''(h_2) \cdot h_{2,xx} \cdot W_3^T) W_2(k,:)^T) \\
		& \;\;\;\; + \sigma'''(h_1(:,k)) \cdot h_{1,x}^2(:,k) \cdot ((\sigma'(h_2) \cdot W_3^T) W_2(k,:)^T) \\
		& \;\;\;\; + \sigma''(h_1(:,k)) \cdot h_{1,xx}(:,k) \cdot ((\sigma'(h_2) \cdot W_3^T) W_2(k,:)^T) .
	\end{split}
\end{equation}

\newpage
Finally, the derivatives with respect to the weights $w_{j,k} \in W_1$, for $j = 1, \dots, m_1$ and $k=1,\dots,m_2$, are given by
\begin{equation}
	\begin{split}
		\frac{\partial \tilde{u}}{\partial w_{j,k}} 
		& = X(:,j) \cdot \sigma'(h_1(:,k)) \cdot ((\sigma'(h_2) \cdot W_3^T) W_2(k,:)^T) \\
		\frac{\partial \tilde{u}_t}{\partial w_{j,k}} 
		& = X(:,j) \cdot \sigma'(h_1(:,k)) \cdot ((\sigma''(h_2) \cdot h_{2,t} \cdot W_3^T) W_2(k,:)^T) \\
		& \;\;\;\; + X(:,j) \cdot \sigma''(h_1(:,k)) \cdot h_{1,t}(:,k) \cdot ((\sigma'(h_2) \cdot W_3^T) W_2(k,:)^T) \\
		& \;\;\;\; + X_t(:,j) \cdot \sigma'(h_1(:,k)) \cdot ((\sigma'(h_2) \cdot W_3^T) W_2(k,:)^T) \\
		\frac{\partial \tilde{u}_x}{\partial w_{j,k}} 
		& = X(:,j) \cdot \sigma'(h_1(:,k)) \cdot ((\sigma''(h_2) \cdot h_{2,x} \cdot W_3^T) W_2(k,:)^T) \\
		& \;\;\;\; + X(:,j) \cdot \sigma''(h_1(:,k)) \cdot h_{1,x}(:,k) \cdot ((\sigma'(h_2) \cdot W_3^T) W_2(k,:)^T) \\
		& \;\;\;\; + X_x(:,j) \cdot \sigma'(h_1(:,k)) \cdot ((\sigma'(h_2) \cdot W_3^T) W_2(k,:)^T) \\
		\frac{\partial \tilde{u}_{xx}}{\partial w_{j,k}} 
		& = X(:,j) \cdot \sigma'(h_1(:,k)) \cdot ((\sigma'''(h_2) \cdot h_{2,x}^2 \cdot W_3^T) W_2(k,:)^T) \\
		& \;\;\;\; + X(:,j) \cdot \sigma''(h_1(:,k)) \cdot h_{1,x}(:,k) \cdot ((2 \sigma''(h_2) \cdot h_{2,x} \cdot W_3^T) W_2(k,:)^T) \\
		& \;\;\;\; + X_x(:,j) \cdot \sigma'(h_1(:,k)) \cdot ((2 \sigma''(h_2) \cdot h_{2,x} \cdot W_3^T) W_2(k,:)^T) \\
		& \;\;\;\; + X(:,j) \cdot \sigma'(h_1(:,k)) \cdot ((\sigma''(h_2) \cdot h_{2,xx} \cdot W_3^T) W_2(k,:)^T) \\
		& \;\;\;\; + X(:,j) \cdot \sigma'''(h_1(:,k)) \cdot h_{1,x}^2(:,k) \cdot ((\sigma'(h_2) \cdot W_3^T) W_2(k,:)^T) \\
		& \;\;\;\; + X_x(:,j) \cdot \sigma''(h_1(:,k)) \cdot 2h_{1,x}(:,k) \cdot ((\sigma'(h_2) \cdot W_3^T) W_2(k,:)^T) \\
		& \;\;\;\; + X(:,j) \cdot \sigma''(h_1(:,k)) \cdot h_{1,xx}(:,k) \cdot ((\sigma'(h_2) \cdot W_3^T) W_2(k,:)^T) \\
		& \;\;\;\; + X_{xx}(:,j) \cdot \sigma'(h_1(:,k)) \cdot ((\sigma'(h_2) \cdot W_3^T) W_2(k,:)^T) .
	\end{split}
\end{equation}

\section{Parameters} \label{sec_parameters}

Tabs.\ \ref{tab_parameters} and \ref{tab_problem_parameters} summarize the general optimization settings and problem-specific configurations used in the numerical experiments.

\begin{table}[t]
	\small
	\centering
	\caption{Optimization and training parameters used in the numerical experiments.}
	\label{tab_parameters}
	\begin{tabular}{|l|c|}
		\hline
		Parameter                               &         Value          \\ \hline
		Initial damping parameter $\mu_0$ (LM)  &         $0.01$         \\ \hline
		Damping update strategy (LM)            & Adaptive, factor of 10 \\ \hline
		Maximum iterations (LM)                 &         $4000$         \\ \hline
		Maximum iterations (BFGS, L-BFGS, Adam) &        $20000$         \\ \hline
		Step tolerance (LM, BFGS, L-BFGS)       &       $10^{-12}$       \\ \hline
		History size (L-BFGS)                   &          $10$          \\ \hline
		Initial learning rate $\alpha$ (Adam)   &       $10^{-3}$        \\ \hline
		$\beta_1$ (Adam)                        &         $0.9$          \\ \hline
		$\beta_2$ (Adam)                        &        $0.999$         \\ \hline
		$\epsilon$ (Adam)                       &       $10^{-8}$        \\ \hline
		Learning-rate decay interval (Adam)     &   $1000$ iterations    \\ \hline
		Learning-rate decay factor (Adam)       &         $0.95$         \\ \hline
	\end{tabular}
\end{table}

\begin{table}[t]
	\small
	\centering
	\caption{Problem-specific configurations used in the numerical experiments.}
	\label{tab_problem_parameters}
	\begin{tabular}{|l|c|c|c|c|}
		\hline
		              &         &                          & Collocation & Activation \\
		PDE           & Problem &         Network          &   Points    &  Function  \\ \hline
		Burgers       & Forward &     NN$(2,25,25,1)$      &    10000    &  $\tanh$   \\ \hline
		Schr\"odinger & Forward & $2\times$NN$(3,25,25,1)$ &    10000    &  $\tanh$   \\ \hline
		Allen--Cahn   & Inverse &     NN$(2,25,25,1)$      &    2000     &  $\tanh$   \\ \hline
		Bratu         & Inverse &     NN$(3,25,25,1)$      &    10000    &  $\tanh$   \\ \hline
	\end{tabular}
\end{table}

\end{document}